\documentclass[12pt]{article}
\usepackage{latexsym, amssymb}
\usepackage{mathbbol}

\DeclareMathAlphabet\gothic{U}{euf}{m}{n}

\makeatletter
\def\eqnarray{\stepcounter{equation}\let\@currentlabel=\theequation
\global\@eqnswtrue
\tabskip\@centering\let\\=\@eqncr
$$\halign to \displaywidth\bgroup\hfil\global\@eqcnt\z@
  $\displaystyle\tabskip\z@{##}$&\global\@eqcnt\@ne
  \hfil$\displaystyle{{}##{}}$\hfil
  &\global\@eqcnt\tw@ $\displaystyle{##}$\hfil
  \tabskip\@centering&\llap{##}\tabskip\z@\cr}

\def\endeqnarray{\@@eqncr\egroup
      \global\advance\c@equation\m@ne$$\global\@ignoretrue}
\makeatother

\begin{document}
 \bibliographystyle{tom}

\newtheorem{lemma}{Lemma}[section]
\newtheorem{thm}[lemma]{Theorem}
\newtheorem{cor}[lemma]{Corollary}
\newtheorem{voorb}[lemma]{Example}
\newtheorem{rem}[lemma]{Remark}
\newtheorem{prop}[lemma]{Proposition}
\newtheorem{stat}[lemma]{{\hspace{-5pt}}}

\newenvironment{remarkn}{\begin{rem} \rm}{\end{rem}}
\newenvironment{exam}{\begin{voorb} \rm}{\end{voorb}}

\newcommand{\gota}{\gothic{a}}
\newcommand{\gotb}{\gothic{b}}
\newcommand{\gotc}{\gothic{c}}
\newcommand{\gote}{\gothic{e}}
\newcommand{\gotf}{\gothic{f}}
\newcommand{\gotg}{\gothic{g}}
\newcommand{\gothh}{\gothic{h}}
\newcommand{\gotk}{\gothic{k}}
\newcommand{\gotm}{\gothic{m}}
\newcommand{\gotn}{\gothic{n}}
\newcommand{\gotp}{\gothic{p}}
\newcommand{\gotq}{\gothic{q}}
\newcommand{\gotr}{\gothic{r}}
\newcommand{\gots}{\gothic{s}}
\newcommand{\gotu}{\gothic{u}}
\newcommand{\gotv}{\gothic{v}}
\newcommand{\gotw}{\gothic{w}}
\newcommand{\gotz}{\gothic{z}}
\newcommand{\gotA}{\gothic{A}}
\newcommand{\gotB}{\gothic{B}}
\newcommand{\gotG}{\gothic{G}}
\newcommand{\gotL}{\gothic{L}}
\newcommand{\gotS}{\gothic{S}}
\newcommand{\gotT}{\gothic{T}}

\newcounter{teller}
\renewcommand{\theteller}{\Roman{teller}}
\newenvironment{tabel}{\begin{list}%
{\rm \bf \Roman{teller}.\hfill}{\usecounter{teller} \leftmargin=1.1cm
\labelwidth=1.1cm \labelsep=0cm \parsep=0cm}
                      }{\end{list}}

\newcounter{tellerr}
\renewcommand{\thetellerr}{\roman{tellerr}}
\newenvironment{subtabel}{\begin{list}%
{\rm  \roman{tellerr}.\hfill}{\usecounter{tellerr} \leftmargin=1.1cm
\labelwidth=1.1cm \labelsep=0cm \parsep=0cm}
                         }{\end{list}}

\newcounter{proofstep}
\newcommand{\nextstep}{\refstepcounter{proofstep}\ruimte \par 
          \noindent{\bf Step \theproofstep} \hspace{5pt}}
\newcommand{\firststep}{\setcounter{proofstep}{0}\nextstep}

\newcommand{\Ni}{{\bf N}}
\newcommand{\Ri}{{\bf R}}
\newcommand{\Ci}{{\bf C}}
\newcommand{\Ti}{{\bf T}}
\newcommand{\Zi}{{\bf Z}}
\newcommand{\Fi}{{\bf F}}

\newcommand{\proof}{\mbox{\bf Proof} \hspace{5pt}} 
\newcommand{\remark}{\mbox{\bf Remark} \hspace{5pt}}
\newcommand{\ruimte}{\vskip10.0pt plus 4.0pt minus 6.0pt}

\newcommand{\ad}{{\mathop{\rm ad}}}
\newcommand{\Ad}{{\mathop{\rm Ad}}}
\newcommand{\Aut}{\mathop{\rm Aut}}
\newcommand{\arccot}{\mathop{\rm arccot}}
\newcommand{\diam}{\mathop{\rm diam}}
\newcommand{\divv}{\mathop{\rm div}}
\newcommand{\codim}{\mathop{\rm codim}}
\newcommand{\RRe}{\mathop{\rm Re}}
\newcommand{\IIm}{\mathop{\rm Im}}
\newcommand{\Tr}{{\mathop{\rm Tr}}}
\newcommand{\supp}{\mathop{\rm supp}}
\newcommand{\sgn}{\mathop{\rm sgn}}
\newcommand{\essinf}{\mathop{\rm ess\,inf}}
\newcommand{\esssup}{\mathop{\rm ess\,sup}}
\newcommand{\Int}{\mathop{\rm Int}}
\newcommand{\Leibniz}{\mathop{\rm Leibniz}}
\newcommand{\lcm}{\mathop{\rm lcm}}
\newcommand{\loc}{{\rm loc}}
\newcommand{\rlim}{\mathop{\rm r.lim}}

\newcommand{\mod}{\mathop{\rm mod}}
\newcommand{\spann}{\mathop{\rm span}}
\newcommand{\ubar}{\underline{\;}}
\newcommand{\one}{\mathbb{1}}

\hyphenation{groups}
\hyphenation{unitary}

\newcommand{\cb}{{\cal B}}
\newcommand{\cc}{{\cal C}}
\newcommand{\cd}{{\cal D}}
\newcommand{\ce}{{\cal E}}
\newcommand{\cf}{{\cal F}}
\newcommand{\ch}{{\cal H}}
\newcommand{\ci}{{\cal I}}
\newcommand{\ck}{{\cal K}}
\newcommand{\cl}{{\cal L}}
\newcommand{\cm}{{\cal M}}
\newcommand{\co}{{\cal O}}
\newcommand{\cs}{{\cal S}}
\newcommand{\ct}{{\cal T}}
\newcommand{\cx}{{\cal X}}
\newcommand{\cy}{{\cal Y}}
\newcommand{\cz}{{\cal Z}}

\thispagestyle{empty}

\begin{center}
{\Large\bf Small time asymptotics of diffusion processes  } \\[5mm]
\large A.F.M. ter Elst$^1$, Derek W. Robinson$^2$ and Adam Sikora$^3$

\end{center}

\vspace{5mm}

\begin{center}
{\bf Abstract}
\end{center}

\begin{list}{}{\leftmargin=1.8cm \rightmargin=1.8cm \listparindent=10mm 
   \parsep=0pt}
\item
We establish the short-time asymptotic behaviour of the Markovian semigroups
associated with strongly local Dirichlet forms under very general hypotheses.
Our results apply to a wide class of strongly elliptic, subelliptic and degenerate
elliptic operators.
In the degenerate case the asymptotics incorporate possible non-ergodicity.
\end{list}

\vspace{7cm}
\noindent
November 2005

\vspace{5mm}

\vspace{5mm}

\noindent
{\bf Home institutions:}    \\[3mm]
\begin{tabular}{@{}cl@{\hspace{10mm}}cl}
1. & Department of Mathematics  & 
  2. & Centre for Mathematics   \\
& \hspace{15mm} and Computing Science & 
  & \hspace{15mm} and its Applications  \\
& Eindhoven University of Technology & 
  & Mathematical Sciences Institute  \\
& P.O. Box 513 & 
  & Australian National University  \\
& 5600 MB Eindhoven & 
  & Canberra, ACT 0200  \\
& The Netherlands & 
  & Australia  \\[8mm]
3. & Department of Mathematical Sciences & 
 {}\\
& New Mexico State University & 
  {}\\
& P.O. Box 30001 & 
 {}\\
& Las Cruces & 
 {} \\
& NM 88003-8001, USA & 
 {} \\
& {} & 
  & {} 
\end{tabular}

\newpage
\setcounter{page}{1}

\section{Introduction}\label{Sas1}

One of the iconic results in the theory of second-order elliptic operators is
Varadhan's \cite{Varadhan1} \cite{Varadhan2} identification of the small time asymptotic limit
\begin{equation}
\lim_{t \downarrow 0}\,t\log K_t(x\,;y)=-4^{-1}d(x\,;y)^2
\label{eas1.1}
\end{equation}
of the heat kernel $K$ of a strongly elliptic operator on $\Ri^d$
in terms of the intrinsic Riemannian distance $d(\cdot\,;\cdot)$.
The small time behaviour  was subsequently analyzed 
at length by Molchanov \cite{Mol} who extended  (\ref{eas1.1}) 
to a much wider class of operators and manifolds.
These results were then analyzed, largely by probabilistic methods, by various authors
(see, for example, the Paris lectures \cite{Aze2}).
Most of the early results were restricted  to non-degenerate operators with smooth coefficients
acting on smooth manifolds. 
Optimal results for the heat flow on Lipschitz Riemannian manifolds were  obtained much later by
Norris \cite{Nor1}.
The asymptotic relation (\ref{eas1.1}) has, however, been   established for certain classes of 
degenerate subelliptic operators by several authors, in  particular 
for sublaplacians constructed from vector fields satisfying H\"ormander's 
condition for hypoellipticity \cite{BKRR} \cite{Lea1} \cite{Lea2} \cite{KuS}.
It is nevertheless  evident from explicit examples that (\ref{eas1.1}) fails for large classes
of degenerate elliptic operators.
Difficulties arise, for example, from non-ergodic behaviour.

The problems introduced by degeneracies are illustrated by the operator
 $H=-d\,c_\delta\,d$, where $d=d/dx$ and 
\begin{equation}
c_\delta(x)
=\Big( \frac{x^2}{1+x^2} \Big)^\delta
\;\;\; ,  
\label{eSas1;97}
\end{equation}
acting on the real line.
If $\delta\in[0,1/2\rangle$ then the associated diffusion is ergodic.
If, however, $\delta\geq 1/2$ then the diffusion process separates into two processes
on the left and 
right half lines, respectively (see \cite{ERSZ1}, Proposition~6.5).
The degeneracy of $c_\delta$ at the origin creates an impenetrable obstacle for the diffusion.
Therefore the corresponding kernel $K$  satisfies $K_t(x\,;y)=0$ for all $x<0$, $y>0$ and $t>0$
and the effective distance between the left and right half lines is infinite.
But this behaviour is not reflected by the Riemannian distance $d(x\,;y)=|\int^y_x c_\delta^{-1/2}|$ 
which is finite for all $\delta\in[0,1\rangle$.
Hence (\ref{eas1.1}) must fail for this diffusion process if $\delta\in[1/2,1\rangle$.
More complicated  phenomena  can occur for degenerate operators in higher dimensions.
To be specific let $I = \{ (\alpha,0) : \alpha \in [-1,1] \} $ 
be a bounded one-dimensional interval in $\Ri^2$ and consider the form 
$h(\varphi)=\sum^2_{i=1}(\partial_i\varphi,\,c_{\delta}\,\partial_i\varphi)$ 
with $D(h)=W^{1,2}(\Ri^2)$
where  $c_{\delta}(x)=(|x|^2_I/(1+|x|_I^2))^\delta$
and $|x|_I$  denotes the Euclidean distance from $x\in\Ri^2$ to the interval $I$.
It follows that if $\delta\geq 1/2$ then the interval $I$  
presents an impenetrable obstacle for the corresponding diffusion 
and the effective configuration space for the 
process is $\Ri^2\backslash I$.
In particular the appropriate distance for the description of the diffusion
 is the intrinsic Riemannian distance  on 
$\Ri^2\backslash I$ rather than that on $\Ri^2$.
Although the Riemannian distance on $\Ri^2$ is well-defined if $\delta\in [0,1\rangle$
it is  not suited to the description of the diffusion if $\delta\in[1/2,1\rangle$.
Thus  the problem of a deeper  understanding of the small time behaviour
of the heat kernel associated with degenerate  operators consists in part in 
identifying the appropriate measure of distance.

Hino and Ram\'{\i}rez \cite{HiR} (see also  \cite{Hin1} \cite{Ram}) made considerable progress in understanding the small
time asymptotics of general  diffusion processes by examining the problem in the broader
context of Dirichlet forms on a $\sigma$-finite  measure space $(X,\cb, \mu)$ \cite{FOT} \cite{BH} \cite{Mosco}.
First they consider  an integrated version of (\ref{eas1.1}).
Set $K_t(A\,;B)=\int_Ad\mu(x)\int_Bd\mu(y)\,K_t(x\,;y)$ for measurable subsets $A,B\in\cb$.
Then the pointwise asymptotic estimate (\ref{eas1.1}) leads, under quite general conditions,
to a set-theoretic version
\begin{equation}
\lim_{t \downarrow 0}\,t\log K_t(A\,;B)=-4^{-1}d(A\,;B)^2
\label{eas1.2}
\end{equation}
for open sets where the distance between the sets $A$ and $B$ is defined in the usual manner with infima.
Secondly, if the measure $\mu$ is finite then  Hino--Ram\'{\i}rez establish (\ref{eas1.2}) for the  kernel of the 
semigroup  corresponding to a strongly local, conservative,  Dirichlet form on $L_2(X\,;\mu)$
and for bounded measurable sets $A$ and $B$ but with a
 set-theoretic distance $d(A\,;B)$  defined directly in terms of the Dirichlet form.
This distance  takes  values in $[0,\infty]$ and is not necessarily the distance arising from
any underlying Riemannian structure.
Nevertheless the estimate (\ref{eas1.2}) establishes that 
it is the correct measure of small time behaviour.
A key feature of this formalism is that it allows for the possibility that $K_t(A\,;B)=0$ 
for all small $t$ and then $d(A\,;B)=\infty$.
The principal disadvantages of the Hino--Ram\'{\i}rez result is that it requires
 $(X,\cb, \mu)$ to be  a probability space and the form to be conservative.
One of our aims is to remove these restrictions and to derive the estimate (\ref{eas1.2}) for diffusion
processes related to a large class of regular, strongly local, Markovian forms on a general measure
space $(X,\cb, \mu)$.
Our formalism is suited to applications to degenerate elliptic operators.
A second aim is to establish conditions which allow one to pass from the estimate (\ref{eas1.2}) 
to the pointwise estimate (\ref{eas1.1}).
We prove a general result which covers a variety  of canonical situations.
For example,  we are able to derive
pointwise estimates for subelliptic operators on Lie groups, to give an independent proof of Norris' result \cite{Nor1}
for the Laplace--Beltrami operator on a Lipschitz Riemannian manifold and to discuss Schr\"odinger-like semigroups
with locally bounded potentials.

Simultaneous with this work Ariyoshi and Hino \cite{AH} extended the Hino--Ram\'{\i}rez result to strongly 
local Dirichlet forms on general $\sigma$-finite measure spaces.
Their tactic is broadly similar to ours although both were developed independently.
Both proofs rely on local approximations and limits of the local diffusions and the corresponding distances.
The limits are controlled by use of the evolution equation associated with the form.
Ariyoshi and Hino refine the  estimates on the equations of motion given in \cite{HiR} and \cite{Ram} whilst our
arguments rely  on the observation that the corresponding wave equation has a finite speed of propagation \cite{ERSZ1}.
The latter property gives rather precise information about the small scale evolution.

In Section~\ref{Sas2} we establish the general formalism in which we work and give a 
complete definition of the distance $d(A\,;B)$.
Then we are able to give a precise statement of our conclusions.
Proofs are given in Section~\ref{Sas3}.
In Section~\ref{Sas4} we discuss various applications to elliptic and subelliptic differential
operators.

\section{General formalism} \label{Sas2}

We begin by summarizing some standard definitions and results on Markovian forms and Dirichlet forms.
As background we refer to  the books by Fukushima, Oshima and Takeda \cite{FOT},
Bouleau and Hirsch \cite{BH} and  Ma and R\"ockner \cite{MR} 
and the papers by LeJan \cite{LJ1} and Mosco \cite{Mosco}.
Most of the early results we need are summarized in  Mosco \cite{Mosco}
together with more recent results, notably on relaxed forms and convergence properties.

Let $X$ be a connected locally compact separable
metric space equipped with a positive Radon measure $\mu$, such that  $\supp \mu = X$.
We consider forms and operators on the real Hilbert space $L_2(X)$ and all functions 
in the sequel are assumed to be real-valued.
Moreover, all forms are assumed to be symmetric and densely defined but not necessarily closed or even closable.
We adopt the definition of Markovian form given in \cite{FOT}, page~4, or 
Mosco \cite{Mosco}, page~374.
Then a Dirichlet form
is a closed Markovian form.

A Dirichlet form $\ce$ on $L_2(X)$ is called {\bf regular} if there is a subset of $D(\ce) \cap C_c(X)$
which is  a core of $\ce$, i.e., which is dense  in $D(\ce)$
with respect to the natural norm $\varphi \mapsto (\ce(\varphi) + \|\varphi\|_2^2)^{1/2}$,
and which is also dense  in $C_0(X)$ with respect to the supremum norm.
There are several non-equivalent kinds of locality for forms or Dirichlet forms 
\cite{BH} \cite{FOT} \cite{Mosco}.
In general we will call any positive quadratic form $h$ {\bf strongly local}
if  $h(\psi,\varphi) = 0$ for all $\varphi,\psi \in D(h)$ with $\supp \varphi$ and 
$\supp \psi$ compact and $\psi$  constant on a neighbourhood of $\supp \varphi$.
Following Mosco \cite{Mosco} we refer to a strongly local regular Dirichlet form 
as a {\bf diffusion}.

There is a possibly stronger version of the locality condition given in \cite{BH}, Section~I.5.
A Dirichlet form $\ce$ is called {\bf [BH]-local} if $\ce(\psi,\varphi) = 0$ for all $\varphi,\psi \in D(\ce)$
and $a \in \Ri$ such that $(\varphi + a \one) \psi = 0$.
It follows, however, by \cite{BH}
Remark~I.5.1.5 and Proposition~I.5.1.3 $(L_0) \Rightarrow (L_2)$, that strong locality and
[BH]-locality are equivalent for regular Dirichlet forms.
Therefore diffusions
are [BH]-local.

The Dirichlet form $\ce$ is called {\bf conservative} if $\one \in D(\ce)$ and 
$\ce(\one) = 0$.
The condition $\one \in D(\ce)$ of course requires that $\one \in L_2(X)$, i.e., 
$\mu(X)<\infty$,
so it is of limited applicability.

Let $\ce$ be a  general Dirichlet form on $L_2(X)$.
First, 
for all $\psi \in D(\ce) \cap L_\infty(X)$ define 
$\ci^{(\ce)}_\psi \colon D(\ce) \cap L_\infty(X) \to \Ri$ by
\[
\ci^{(\ce)}_\psi(\varphi) = \ce(\psi \, \varphi,\psi) - 2^{-1} \ce(\psi^2,\varphi)
\;\;\; .  
\]
If no confusion is possible  we drop the suffix and write 
$\ci_\psi(\varphi) = \ci^{(\ce)}_\psi(\varphi)$.
If $\varphi\geq 0$ it follows that 
$\psi \mapsto \ci_\psi(\varphi)$
is a Markovian form with  domain $D(\ce)\cap L_\infty(X)$
which satisfies the key properties
\begin{equation}
0\leq \ci_{F\circ\psi}(\varphi) \leq \ci_\psi(\varphi) \leq \|\varphi\|_\infty \, \ce(\psi)
\label{easr1.11}
\end{equation}
for all $\varphi,\psi \in D(\ce) \cap L_\infty(X)$ with $\varphi\geq0$ and all normal contractions
$F$ (see  \cite{BH}, Proposition~I.4.1.1).
Hence
\begin{equation}
|\ci_{\psi_1}(\varphi) - \ci_{\psi_2}(\varphi)|
\leq \ci_{\psi_1 - \psi_2}(\varphi)^{1/2} \, \ci_{\psi_1 + \psi_2}(\varphi)^{1/2}
\leq \ce(\psi_1 - \psi_2)^{1/2} \, \ce(\psi_1 + \psi_2)^{1/2} \, \|\varphi\|_\infty
\label{easr1.115}
\end{equation}
for all $\varphi,\psi_1,\psi_2 \in D(\ce) \cap L_\infty(X)$ with $\varphi\geq0$.
This form  is referred to as the truncated form by Roth
 \cite{Roth}, Theorem~5.

A regular Dirichlet form has a canonical representation originating with Beurling and
Deny \cite{BeD}  and in its final form by LeJan \cite{LJ1} 
(see \cite{Mosco}, Section~3e or \cite{FOT}, Section~3.2). 
If $\ce$ is strongly local, i.e., if it is a diffusion, this representation takes the simple form
\[
\ce(\psi) =\int_X d\mu_\psi(x)
\]
for all $\psi\in D(\ce)$ where the $\mu_\psi\,(=\mu^{(\ce)}_\psi\,)$ are positive Radon measures.
The measures are uniquely determined by the identity
\[
\ci_\psi(\varphi) =\int_X d\mu_\psi \, \varphi
\]
for all $\varphi, \psi \in D(\ce) \cap L_\infty(X)$
and the continuity property
\[
\|\mu_{\psi_1} - \mu_{\psi_2}\|
\leq \ce(\psi_1 - \psi_2)^{1/2} \, \ce(\psi_1 + \psi_2)^{1/2}
\]
for all $\psi_1,\psi_2 \in D(\ce)$ as a consequence of (\ref{easr1.115})
(see \cite{LJ1}, Propositions~1.4.1 and 1.5.1 or \cite{FOT},
Section~3.2 and Lemma~3.2.3).
The measure $\mu_\psi$ is usually referred to as the {\bf energy measure}.

If $\ce$ and $\cf$ are two Dirichlet forms which are [BH]-local and satisfy $\ce\leq \cf$ then
\begin{equation}
\ci^{(\ce)}_\psi(\varphi)\leq \ci^{(\cf)}_\psi(\varphi)
\label{easr1.12}
\end{equation} 
for all $\varphi,\psi\in D(\cf)\cap L_\infty(X)$ with $\varphi\geq0$
(see, for example, \cite{ERSZ2} Proposition~3.2).

\smallskip

Let $\ce$ be a diffusion. 
Define $D(\ce)_{\rm loc}$ as the vector space of 
(equivalent classes of) all measurable functions $\psi \colon X \to \Ci$
such that for every compact subset $K$ of $X$ there is a $\hat \psi \in D(\ce)$
with $\psi|_K = \hat \psi|_K$.
Since $\ce$ is   regular and [BH]-local one can define 
$\widehat \ci_\psi^{(\ce)} = \widehat \ci_\psi \colon D(\ce)\cap L_{\infty,c}(X) \to \Ri$ by
\[
\widehat \ci_\psi(\varphi) = \ci_{\hat \psi}(\varphi)
\]
for all $\psi \in D(\ce)_{\rm loc} \cap L_\infty(X)$ and 
$\varphi \in D(\ce) \cap L_{\infty,c}(X)$ where 
$\hat \psi \in D(\ce)\cap L_\infty(X)$ is such that 
$\psi|_{\supp \varphi} = \hat \psi|_{\supp \varphi}$.
Next, for all $\psi \in D(\ce)_{\rm loc} \cap L_\infty(X)$, define
\[
|||\widehat \ci_\psi|||
= \sup \{ \,|\widehat \ci_\psi(\varphi)| : 
      \varphi \in D(\ce) \cap L_{\infty,c}(X) , \; \|\varphi\|_1 \leq 1\, \}
\in [0,\infty]
\;\;\; .  \]
Now, for all $\psi \in L_\infty(X)$ and measurable sets $A,B \subset X$, introduce 
\begin{eqnarray*}
d_\psi(A\,;B)
&=& \sup \{ \,M \in \Ri : \psi(x) - \psi(y) \geq M \mbox{ for a.e.\ $x \in A$ and a.e.\ }
                      y \in B\, \}  \\[5pt]
&=&\essinf_{x\in A}\psi(x)-\esssup_{y\in B}\psi(y) 
\in \langle -\infty,\infty]
\;\;\; .  
\end{eqnarray*}
Recall that 
\begin{eqnarray*}
\esssup_{y\in B}\psi(y) 
& = & \inf \{\, m \in \Ri : | \{ y \in B : \psi(y) > m \} | = 0 \,\}   \\[0pt]
& = & \min \{\, m \in [-\infty,\infty \rangle : | \{ y \in B : \psi(y) > m \,\} | = 0 \} 
\in [-\infty,\infty \rangle
\end{eqnarray*}
and $\essinf_{x \in A} \psi(x) = - \esssup_{x \in A} -\psi(x)$.
Finally define
\[
d(A\,;B)
= d^{(\ce)}(A\,;B)
= \sup \{ \,d_\psi(A\,;B) : \psi \in D_0(\ce) \, \}
\;\;\;,
\]
where
\[
D_0(\ce)
= \{ \psi\in D(\ce)_{\rm loc} \cap L_\infty(X): |||\widehat \ci_\psi||| \leq 1 \}
\]
as in \cite{ERSZ2}.
A similar definition was given by Hino and Ram\'{\i}rez \cite{HiR} (see also \cite{Stu2}),
but since they considered probability spaces the introduction of $D(\ce)_{\rm loc}$ was 
unnecessary.
Regularity of the form was also unnecessary.
If, however, we were to replace $D(\ce)_{\rm loc}$ by $D(\ce)$ in the definition 
of $d(A\,;B)$ then, since $\psi\in L_2(X)$, one would obtain $d_\psi(A\,;B)\leq  0$ for all measurable sets $A,B \subset X$
with $|A| = |B| = \infty$ and this would give $d(A\,;B)= 0$.
On the other hand  the definition with $D(\ce)_{\rm loc}$ is not useful unless
$D(\ce)$ contains sufficient bounded functions with compact support.
Hence regularity of some sort is essential.

The distance determines the following $L_2$-off diagonal bounds, or Davies--Gaffney bounds
on the semigroup associated to a diffusion by Theorem~1.2 in \cite{ERSZ2}.

\begin{thm} \label{tdfd3.1}
Let $\ce$ be a diffusion over $X$. 
If  $S^{(\ce)}$ denotes the semigroup generated by the self-adjoint  operator $H$
on $L_2(X)$ associated with $\ce$ and if 
 $A$ and $B$ are measurable subsets of~$X$ then
\[
|(\psi, S^{(\ce)}_t\varphi)|\leq e^{-d^{(\ce)}(A;B)^2(4t)^{-1}}\|\psi\|_2 \, \|\varphi\|_2
\]
for all $\psi\in L_2(A)$, $\varphi \in L_2(B)$ and $t>0$ with the convention
$e^{-\infty}=0$.
\end{thm}

Let  $h$ be a positive quadratic form on $L_2(X)$.
Then  $h$ is not necessarily closable but there exists a
largest closed quadratic form,  denoted by $h_0$,  which is majorized by $h$, i.e.,
$D(h) \subseteq D(h_0)$ and $h_0(\varphi) \leq h(\varphi)$
for all $\varphi \in D(h)$.
The form $h_0$ is referred to in the literature on discontinuous media as the {\bf relaxation}
of $h$ (see \cite{Mosco}, Section~1.e).
It can be understood in various ways.
Simon (see \cite{bSim5}, Theorems~2.1 and 2.2) established    that 
 $h$ can be decomposed uniquely as a sum $h=h_r+h_s$ of two positive forms with 
$D(h_r)=D(h)=D(h_s)$ and
$h_r$ the largest closable form majorized by~$h$.
Simon refers to $h_r$ as the {\bf regular part} of $h$.
Then $h_0=\overline{h_r}$, the closure of $h_r$.
Simon also proved that $D(h_0)$ consists of those $\varphi\in L_2(X)$
for which there is a sequence $\varphi_n\in D(h)$ such that 
$\lim_{n \to \infty} \varphi_n = \varphi$
in $L_2(X)$ and $\liminf_{n\to\infty} h(\varphi_n)<\infty$.
Moreover, $h_0(\varphi)$ equals the minimum of all $\liminf_{n \to \infty} h(\varphi_n)$, 
where the minimum is taken over all $\varphi_1,\varphi_2,\ldots \in D(h)$ such that 
$\lim_{n \to \infty} \varphi_n = \varphi$ in $L_2(X)$.
(See \cite{bSim4}, Theorems~2 and~3.)
Note that if $D$ is a subspace of $L_2(X)$ which is dense in $D(h)$ then 
$(h|_D)_0 = h_0$.
Moreover, if $h$ and $k$ are two positive quadratic forms with $h \leq k$ then $h_0 \leq k_0$.

The relaxation can also be understood by approximation.
If $h,h_1,h_2,\ldots$ are closed positive quadratic forms on $L_2(X)$
and $H, H_1,H_2,\ldots$ 
the corresponding positive self-adjoint operators
then we write $h=\rlim_{n \to \infty} h_n $ if 
 $H_1,H_2,\ldots$ 
 converges to $H$  in the strong resolvent sense, i.e., if
$(I + H)^{-1}=\lim_{n \to \infty} (I + H_n)^{-1} $ strongly.
It follows from \cite{Mosco}, Theorem 2.4.1, that $h=\rlim_{n \to \infty} h_n $
if and only if $h(\varphi) \leq \liminf_{n \to \infty} h_n(\varphi_n)$
for all $\varphi,\varphi_1,\varphi_2,\ldots \in L_2(X)$ with $\lim_{n \to \infty} \varphi_n = \varphi$
weakly in $L_2(X)$ and, in addition, for all $\varphi \in L_2(X)$ there exists
a sequence $\varphi_1,\varphi_2,\ldots \in L_2(X)$ such that $\lim_{n \to \infty} \varphi_n = \varphi$
strongly in $L_2(X)$ and $h(\varphi) = \liminf_{n \to \infty} h_n(\varphi_n)$.

If $h_1,h_2,\ldots$ are closed positive quadratic forms on $L_2(X)$
such that $h_1 \geq h_2 \geq \ldots$ then it follows from a result of Kato \cite{Kat1},
Theorem~VIII.3.11, that the corresponding sequence $H_1,H_2,\ldots$ of operators
converges in the strong resolvent sense to a positive self-adjoint operator $H$.
If $h$ is the form corresponding to $H$, then  $h=\rlim_{n\to\infty}h_n$ and 
 $h$ is the largest closed form which is majorized by $h_n$ for all $n \in \Ni$
(Simon \cite{bSim4}, Theorem~3.2).

Now let  $h$ be  a positive quadratic form on $L_2(X)$ and $l$ a closed positive
quadratic form  such that $h\leq\lambda\,l$, for some
$\lambda >0$, and with $D(l)$  dense in $D(h)$. 
Then the forms $h_\varepsilon=h+\varepsilon\,l$, 
with $\varepsilon>0$, are all closed positive forms, with domain $D(l)$, since 
$\varepsilon\,l\leq h_\varepsilon\leq (\lambda+\varepsilon)\,l$.
But $\varepsilon\mapsto h_\varepsilon$ is monotonically decreasing if $\varepsilon \downarrow 0$
and it follows
from the results of Kato and Mosco cited in the foregoing paragraphs that
$\rlim_{\varepsilon\downarrow 0}h_\varepsilon=h_0$ where $h_0$ is the relaxation of $h$.
This characterization justifies the notation $h_0$.
Note that $h_0$ is independent of the particular $l$ used in this construction
which is akin to the viscosity method of partial differential equations.
So $h_0$ could well be called the {\bf viscosity closure} of $h$.
If $h$ is the form of a pure second-order elliptic operator in divergence
form on $\Ri^d$ and   $l$ is  the form of the Laplacian then
the condition $h\leq\lambda\,l$  corresponds to uniform boundedness of the coefficients.

\smallskip

Throughout the remainder of this paper we fix a  strongly local regular 
Dirichlet form $l$ on $L_2(X)$, i.e., a diffusion, satisfying the following property.

\bigskip

\noindent{\bf Condition L}
{\it 
The function $d^{(e)} \colon X \times X \to [0,\infty]$ defined by
\[
d^{(e)}(x\,;y)
= \sup \{ |\psi(x) - \psi(y)| : \psi \in D(l)_{\rm loc} \cap C_{\rm b}(X) \mbox{ and } |||\widehat \ci^{(l)}_\psi||| \leq 1 \}
  \]
is a metric, the topology induced by this metric 
equals the original topology on $X$ and  the 
balls $B^{(e)}(x\,;r)$ defined by the metric $d^{(e)}$ are relatively compact
for all $x \in X$ and $r > 0$.
}

\ruimte

The first result  of this paper is given by  the following.

\begin{thm} \label{tlbs121}
Let $l$ be a diffusion on $L_2(X)$ 

satisfying {\rm Condition~L}.
Further let $h$ be a positive  form over $X$  such that
$h \leq \lambda \, l_{D(l) \cap L_\infty(X)}$ for some $\lambda > 0$ and  
$h|_{D(l) \cap L_\infty(X)}$ is strongly local and Markovian.

Then the relaxation $h_0$ of $h$ is a diffusion.
Moreover, if $A,B \subset X$ are relatively compact and measurable
then 
\begin{equation}
\lim_{t \downarrow  0} \,t \log (\one_A, S^{(h_0)}_t \one_B) = - 4^{-1} d^{(h_0)}(A\,;B)^2  
\label{etlbs113;1}
\end{equation}
where $S^{(h_0)}$ denotes the semigroup associated with $h_0$.
\end{thm}

Note that if $S^{(h_0)}$ has an integrable kernel $K^{(h_0)}$ then
\[
(\one_A, S^{(h_0)}_t \one_B)=\int_Ad\mu(x)\int_Bd\mu(y)\,K^{(h_0)}_t(x\,;y)=K^{(h_0)}_t(A\,;B)
\]
so (\ref{etlbs113;1}) corresponds to the integrated version (\ref{eas1.2}) of Varadhan's result
(\ref{eas1.1}).

Next note that it follows by choosing $h$ equal to $l$ that one has the following conclusion.

\begin{cor}\label{cas2.1}
If $l$ is a diffusion on $L_2(X)$ 
satisfying {\rm Condition~L}, 
then
\[
\lim_{t \downarrow  0} \,t \log (\one_A, S^{(l)}_t \one_B) = - 4^{-1} d^{(l)}(A\,;B)^2  
\]
for all relatively compact measurable $A,B\subset X$.
\end{cor}

One application of the corollary is to strongly elliptic operators in divergence form on $\Ri^d$
with bounded coefficients. 
Then $d^{(e)}(\cdot\,;\cdot)$ is the corresponding Riemannian distance
and $d^{(l)}(A\,;B)=d^{(e)}(A\,;B)$  if $A$ and $B$ are open and  non-empty,
where  we set
\[
d^{(e)}(A\,;B)=\inf_{x\in A}\inf_{y\in B}d^{(e)}(x\,;y)
\]
for general non-empty $A,B \subseteq X$.

\smallskip

There is also a weaker version of Theorem~\ref{tlbs121} for sets which are possibly not
relatively compact.
In the formulation of this result we let $P_A$ denote the orthogonal projection from 
$L_2(X)$ onto $L_2(A)$.

\begin{cor} \label{cas135}
Assume the conditions of Theorem~{\rm \ref{tlbs121}}.
Then
\[
\lim_{t \downarrow  0} \, t \log \|P_A S^{(h_0)}_t P_B\|_{2 \to 2} 
= - 4^{-1} d^{(h_0)}(A\,;B)^2  
\]
for all measurable $A,B \subset X$.
\end{cor}

It is possible to transform the set-theoretic bounds on the semigroup into pointwise bounds
but this requires some additional assumptions which are satisfied for large classes of
non-degenerate elliptic operators (see Section~\ref{Sas4}).
In general it is possible that $\psi \in D(l)_{\rm loc} \cap L_\infty$ with 
$|||\widehat \ci_\psi||| < \infty$ but $\psi$ is not continuous. 
Therefore one can have  $d^{(l)}(A\,;B)>d^{(e)}(A\,;B)$  even for non-empty open sets $A$ and $B$.
This behaviour is illustrated by the one-dimensional example given in the introduction
 with  $\delta \in [1/2,1\rangle$ in (\ref{eSas1;97}).

This explains the origin of the first assumption in the next theorem.
The second assumption is typically a consequence of a parabolic version of the Harnack inequality.

\begin{thm} \label{tlbs122}
Let $l$ be a diffusion on $L_2(X)$ 
satisfying {\rm Condition~L}.
Assume 
\begin{tabel}
\item\label{tlbs122-1}
 $D_0(l) \subseteq C(X)$, and,  
\item\label{tlbs122-2}
the semigroup $S^{(l)}$ has a continuous kernel~$K$ and 
 there are $\nu,\omega,T \in \langle0,\infty\rangle$ such that 
\[
K_s(x\,;y)
\leq K_t(x\,;z) \, (t s^{-1})^\nu \, e^{\omega(1 + d^{(e)}(y;z)^2 (t-s)^{-1})}
\]
for all $0 < s < t \leq T$ and $x,y,z \in X$.
\end{tabel}
Then 
\[
\lim_{t \downarrow 0} \,t \log K_t(x\,;y) = - 4^{-1} d^{(e)}(x\,;y)^2
\]
for all $x,y \in X$.
\end{thm}

Applications of Theorem~\ref{tlbs122} are discussed in Section~\ref{Sas4}.

\section{Proof of the theorems} \label{Sas3}

First we derive several useful results for general diffusions.

\begin{lemma} \label{llbs110}
Let $\ce$ be a diffusion on $L_2(X)$.
If $A,B$ are measurable with $|A|, |B| < \infty$ and 
\[
\lim_{t \downarrow  0} t \log (\one_A, S^{(\ce)}_t \one_B) = - 4^{-1} d^{(\ce)}(A\,;B)^2
\]
then  $d^{(\ce)}(A\,;B)$ is the supremum of all 
$ r \geq 0$ for which there are $M,t_0 > 0$ such that
\[
(\one_A, S^{(\ce)}_t \one_B) \leq M e^{-(4t)^{-1} r^2} 
\]
 for all $t \in \langle0,t_0]$.
\end{lemma}
\proof\
Let $s$ denote the supremum.
Then the Davies--Gaffney bounds of  Theorem~\ref{tdfd3.1} give $d(A\,;B) \leq s$.
If $d(A\,;B) < s$ then there are $r \in \langle d(A\,;B),s\rangle$ and $M,t_0 > 0$ such that 
$(\one_A, S_t \one_B) \leq M e^{-(4t)^{-1} r^2}$ for all $t \in \langle0,t_0]$.
Then 
\[
\lim_{t \downarrow  0} t \log (\one_A, S_t \one_B) \leq - 4^{-1} r^2 < - 4^{-1} d(A\,;B)^2
\]
which gives a contradiction.\hfill$\Box$

\begin{lemma} \label{las343}
Let $\ce$ be a diffusion on $L_2(X)$.
Let $\psi_1,\psi_2 \in D(\ce)_\loc \cap L_\infty$.
\begin{tabel}
\item  \label{las343-1}
If $\varphi \in D(\ce) \cap L_{\infty,c}$ with $\varphi \geq 0$ then
$\widehat \ci^{(\ce)}_{\psi_1 \wedge \psi_2}(\varphi)
\leq \widehat \ci^{(\ce)}_{\psi_1}(\varphi) + \widehat \ci^{(\ce)}_{\psi_1}(\varphi)$.
\item  \label{las343-2}
$|||\widehat \ci^{(\ce)}_{\psi_1 \vee \psi_2}||| \leq |||\widehat \ci^{(\ce)}_{\psi_1}||| \vee |||\widehat \ci^{(\ce)}_{\psi_2}|||$.
\end{tabel}
\end{lemma}
\proof\
Let $\mu_\psi$ be the energy measure associated with $\ce$.
It follows from~(2.10) in \cite{BM} that 
\[
\mu_{\psi_1 \wedge \psi_2}
= \one_{ \{ x \in X : \psi_1(x) < \psi_2(x) \} } \, \mu_{\psi_1}
   + \one_{ \{ x \in X : \psi_1(x) > \psi_2(x) \} } \, \mu_{\psi_2}
\]
for all $\psi_1,\psi_2 \in D(l) \cap L_\infty$.
Hence 
\begin{eqnarray*}
\ci_{\psi_1 \wedge \psi_2}(\varphi)
& = & \int d\mu_{\psi_1 \wedge \psi_2}\, \varphi 
= \int  d\mu_{\psi_1}\,\varphi \, \one_{ \{ x \in X : \psi_1(x) < \psi_2(x) \} } 
   + \int  d\mu_{\psi_2}\, \varphi \, \one_{ \{ x \in X : \psi_1(x) > \psi_2(x) \} }   \\[5pt]
& \leq & \int d\mu_{\psi_1}\,\varphi 
   + \int  d\mu_{\psi_2}\,\varphi 
= \ci_{\psi_1}(\varphi) + \ci_{\psi_1}(\varphi)
\end{eqnarray*}
for all $\varphi \in D(\ce) \cap L_\infty$ with $\varphi \geq 0$.
Now Statement~\ref{las343-1} of the lemma follows.

The proof of Statement~\ref{las343-2} is similar.\hfill$\Box$

\begin{lemma} \label{las344}
Let $\ce$ be a diffusion on $L_2(X)$.
Let $\psi_1,\psi_2,\ldots \colon X \to [0,\infty\rangle$ be measurable 
and assume that $\psi_n \wedge N \in D_0(\ce)$ for all $n,N \in \Ni$.
Define $\psi \colon X \to [0,\infty]$ by $\psi = \sup_{n \in \Ni} \psi_n$.
Then $\psi \wedge N \in D_0(\ce)$ for all $N \in \Ni$.
\end{lemma}
\proof\
Since $\varphi_1 \vee \varphi_2 \in D_0(\ce)$ for all $\varphi_1,\varphi_2 \in D_0(\ce)$
by Lemma~\ref{las343}.\ref{las343-2} we may assume that 
$\psi_1 \leq \psi_2 \leq \ldots$.
Let $N \in \Ni$.
We may also assume that $\psi_n \leq N$ for all $n \in \Ni$
by (\ref{easr1.11}).

Let $K \subset X$ compact.
Since $\ce$ is regular there exist $\chi,\tilde \chi \in D(\ce) \cap C_c(X)$ 
such that $0 \leq \chi \leq \tilde \chi \leq 1$, $\chi|_K = 1$ and 
$\tilde \chi|_{\supp \chi} = 1$.
It follows from Lemma~\ref{las343}.\ref{las343-1} and the strong locality of $\ce$ that 
\[
\ce(\psi_n \wedge N \chi)
= \ci_{\psi_n \wedge N \chi}(\tilde \chi)  \\[5pt]
\leq \widehat \ci_{\psi_n}(\tilde \chi) + \ci_{N \chi}(\tilde \chi)  \\[5pt]
\leq \|\tilde \chi\|_1 + N^2 \, \ce(\chi)
\]
uniformly for all $n \in \Ni$.
Moreover, $\|\psi_n \wedge N \chi\|_2 \leq N \, \|\chi\|_2$ 
uniformly for all $n \in \Ni$.
Hence the sequence $\psi_1 \wedge N \chi, \psi_2 \wedge N \chi, \ldots$
is bounded in $D(\ce)$.
Therefore there exists a subsequence such that 
$\lim_{k \to \infty} \psi_{n_k} \wedge N \chi$ exists weakly in $D(\ce)$
and $\lim_{k \to \infty} \psi_{n_k} \wedge N \chi$ exists almost everywhere.
By definition of $\psi$ one has 
$\lim_{k \to \infty} \psi_{n_k} \wedge N \chi = \psi \wedge N \chi$
almost everywhere.
Hence $\psi \wedge N \chi \in D(\ce)$.
In particular $\psi \in D(\ce)_\loc \cap L_\infty$.
Moreover, it follows from Statement (a) on page 269 of \cite{bSim4} that 
\begin{eqnarray*}
\widehat \ci_\psi(\varphi)
 =  \ci_{\psi \wedge N \chi}(\varphi)
&\leq& \liminf_{k \to \infty} \ci_{\psi_{n_k} \wedge N \chi}(\varphi)\\[5pt]
&=& \liminf_{k \to \infty} \widehat \ci_{\psi_{n_k}}(\varphi)
\leq \liminf_{k \to \infty} |||\widehat \ci_{\psi_{n_k}}||| \, \|\varphi\|_1
\leq \|\varphi\|_1
\end{eqnarray*}
for all $\varphi \in D(\ce) \cap L_{\infty,c}$ with $\supp \varphi \subset K$ 
and $\varphi \geq 0$.
So $\psi \in D_0(\ce)$.\hfill$\Box$

\begin{lemma} \label{ldfd3n11}
Let $\ce$ be a diffusion on $L_2(X)$.
Let $X_1 \subset X_2 \subset \ldots$ be measurable subsets of $X$ such that $X = \bigcup_{n=1}^\infty X_n$.
If $A,B$ are measurable subsets of $X$ then
\[
d^{(\ce)}(A\,;B) = \lim_{n \to \infty} d^{(\ce)}(A \cap X_n\,;B) = \inf_{n \in \Ni} d^{(\ce)}(A \cap X_n\,;B)
\;\;\; .  \]
\end{lemma}
\proof\
Obviously $d(A \cap X_1\,;B) \geq d(A \cap X_2\,;B) \geq \ldots \geq d(A\,;B)$.
Hence $\lim_{n \to \infty} d(A \cap X_n\,;B) = \inf_{n \to \infty} d(A \cap X_n\,;B) \geq d(A\,;B)$.
Let $M \in [0,\infty\rangle$ and suppose that 
$M \leq \inf_{n \in \Ni} d(A \cap X_n\,;B)$.

Let $\varepsilon > 0$.
Then by Lemma~4.2.III in \cite{ERSZ2}  for all $n \in \Ni$ there exists a $\psi_n \in D_0(\ce)$ 
such that $0 \leq \psi_n \leq M$,
$\psi_n|_B = 0$  and $\psi_n|_{A \cap X_n} \geq M - \varepsilon$.
Set $\psi = \sup_{n \in \Ni} \psi_n$. 
Then $\psi \in D_0(\ce)$ by Lemma~\ref{las344}.
Moreover, $\psi|_A \geq M-\varepsilon$ and $\psi|_B = 0$.
So
\[
M-\varepsilon \leq d_\psi(A\,;B) \leq d(A\,;B)
\;\;\; .  \]
Since this is valid for all $M \in [0,\infty\rangle$ with
$M \leq \inf_{n \in \Ni} d(A \cap X_n\,;B)$ and all $\varepsilon > 0$
it follows that $d(A\,;B) \geq \inf_{n \in \Ni} d(A \cap X_n\,;B)$
and the proof of 
the lemma is complete.\hfill$\Box$

\ruimte

Next we develop some structural results which depend on the domination property assumed in  Theorem~\ref{tlbs121}.
Initially we do not need to assume that $l$ satisfies Condition~L.

We fix from now on a diffusion $l$ over $X$ and set $\gotA=D(l) \cap L_\infty(X)$.
Define   $\cc_l$ as the cone of all positive quadratic forms $h$ over $X$ with $\gotA \subseteq D(h)$  
satisfying  
\begin{enumerate}
\item
$ h \leq \lambda \, l|_\gotA$ for some $\lambda > 0$,
\item
$\gotA$ is dense in $D(h)$, and,
\item
$h|_\gotA$ is strongly local and Markovian.
\end{enumerate}

Note that we do not assume that the forms $h$ are closed or even closable.

The first lemma shows that there is no confusion possible for $h(\varphi)$ 
if $\varphi \in D(l) \backslash \gotA$ under a mild condition that is satisfied if 
$h \in \cc_l$.

\begin{lemma} \label{las341}
Let $h$ be a positive quadratic form over $X$ with $D(h) = \gotA$ and suppose there exists a
$\lambda > 0$ such that $ h \leq \lambda \, l|_\gotA$.
Then there exists a unique positive quadratic form $\tilde h$ over $X$ with $D(\tilde h) = D(l)$ 
and $\tilde h|_\gotA = h$.
This form $\tilde h$ then satisfies $\tilde h \leq \lambda \, l$.
\end{lemma}
\proof\
The existence and the second part of the lemma are easy.

So it remains to prove the uniqueness.
Let $\tilde h,\hat h$ be two positive quadratic form over $X$ with $D(\tilde h) = D(l) = D(\hat h)$ 
and $\tilde h|_\gotA = h = \hat h|_\gotA$.
Let $\varepsilon > 0$.
Set
\[
\tilde h_\varepsilon=(1+\varepsilon)^{-1} (\tilde h + \varepsilon \lambda \, l)
\;\;\;\; \mbox{and} \;\;\;\;
\hat h_\varepsilon=(1+\varepsilon)^{-1} (\hat h + \varepsilon \lambda \, l)
\]
Then $D(\tilde h_\varepsilon) = D(l) = D(\hat h_\varepsilon)$ and 
$\tilde h_\varepsilon(\varphi) = \hat h_\varepsilon(\varphi)$ for all $\varphi \in \gotA$.
Moreover, $\tilde h_\varepsilon(\varphi) \leq \lambda \, l(\varphi)$ for all $\varphi \in \gotA$.
Let $\varphi \in D(l)$.
Since $\gotA$ is dense in $D(l)$ there are $\varphi_1,\varphi_2,\ldots \in \gotA$ such that
$\lim_{n \to \infty} (l(\varphi - \varphi_n) + \|\varphi - \varphi_n\|_2) = 0$.
Since $\tilde h_\varepsilon$ is closed and $\tilde h|_\gotA \leq l|_\gotA$ 
it follows from \cite{Kat1} Theorem~VI.1.12,
that $\tilde h_\varepsilon(\varphi) = \lim_{n \to \infty} \tilde h_\varepsilon(\varphi_n)$.
But similarly $\hat h_\varepsilon(\varphi) = \lim_{n \to \infty} \hat h_\varepsilon(\varphi_n)$.
Therefore $\tilde h_\varepsilon(\varphi) = \hat h_\varepsilon(\varphi)$.
This is valid for all $\varepsilon > 0$.
Hence $\tilde h(\varphi) = \lim_{\varepsilon \downarrow 0} \tilde h_\varepsilon(\varphi) = \ldots = \hat h_\varepsilon(\varphi)$.\hfill$\Box$

\ruimte

\begin{cor} \label{cas343}
If $h \in \cc_l$ then there is a unique $\tilde h \in \cc_l$ such that 
$D(\tilde h) = D(l)$ and $\tilde h|_\gotA = h|_\gotA$.
\end{cor}

It also follows from Lemma~\ref{las341} that for a positive quadratic form $h$ over $X$ with $D(h) \subseteq D(l)$ one can give 
a characterization of the condition $h \in \cc_l$.

\begin{cor} \label{las342}
If $h$ is a quadratic form over $X$ with $D(l) \subseteq D(h)$ 
and if $\lambda>0$ then the 
following are equivalent.
\begin{tabel}
\item
$ h \leq \lambda \, l|_\gotA$ and 
$\gotA$ is dense in $D(h)$.

\item
 $ h \leq \lambda \, l$ and 
\ $D(l)$ is dense in $D(h)$.
\end{tabel}
\end{cor}

Next we begin the analysis of the relaxations of the forms in $\cc_l$.

\begin{lemma} \label{llbs231}
Let $h\in \cc_l$. 
Then the relaxation  $h_0$ of $h$ is a regular Dirichlet form with 
$D(l) \subseteq D(h_0)$.
Moreover, every core for the form~$l$ is a core for $h_0$.
\end{lemma}
\proof\
By assumption the algebra $\gotA$ is dense in $D(h)$.
Therefore $(h|_\gotA)_0 = h_0$.
Since  $h|_\gotA$ is Markovian it follows from 
\cite{Mosco}, Corollary~2.8.2, that $h_0$ is a Dirichlet form.

Let $h_{\gotA \, r}$ be the regular part of $h|_\gotA$ as in \cite{bSim5}.
Then $\gotA = D(h_{\gotA \, r})$.
Hence $\gotA$ is dense in  $D(\overline{h_{\gotA \, r}})$.
But $h_0 = (h|_\gotA)_0$ is the closure $\overline{h_{\gotA \, r}}$ of the regular part $h_{\gotA \, r}$ of $h|_\gotA$.
So $\gotA$ is a core for $h_0$.

Since $h \leq \lambda \, l|_\gotA$ it follows that $h_0 \leq \lambda \, (l|_\gotA)_0 = \lambda \, l$.
But by the above $D(l)$ is a core for~$h_0$.
It then easily follows that any core for $l$ is also a core for~$h_0$.

Finally, since $D(l) \cap C_c(X)$ is dense in $C_0(X)$
and $D(l) \cap C_c(X) \subset D(h_0) \cap C_c(X)$ it follows that 
$D(h_0) \cap C_c(X)$ is dense both in $D(h_0)$ and in $C_0(X)$.\hfill$\Box$

\ruimte

Next set $h_\varepsilon= \tilde h + \varepsilon \, l$ for all $\varepsilon>0$ and $h \in \cc_l$,
where $\tilde h$ is the unique form given by Corollary~\ref{cas343}.
Then $D(h_\varepsilon) = D(l)$.

\begin{cor} \label{clbs231.5}
If $h\in \cc_l$ then $h_\varepsilon$ is a diffusion
for all $\varepsilon > 0$.
Moreover, $\rlim_{\varepsilon \downarrow 0} h_\varepsilon = h_0$.
\end{cor}
\proof\
If $ h \leq \lambda \, l|_\gotA$ then
$\varepsilon \, l \leq h_\varepsilon \leq (\lambda+\varepsilon)\,l$.
Since   $l$ is closed
it follows that $h_\varepsilon$ is closed.
Therefore  $h_\varepsilon = (h_{\varepsilon})_0$
is a regular Dirichlet form by Lemma~\ref{llbs231}.
Moreover, $\gotA$ is a core for $h_\varepsilon$ and 
$h_\varepsilon|_\gotA$ is strongly local.
Therefore $h_\varepsilon$ is strongly local by \cite{BH},
Remark~I.5.1.5 and Proposition~I.5.1.3 $(L_0) \Rightarrow (L_2)$.

If $\tilde h$ is the unique form as in Corollary~\ref{cas343} then
\[
\rlim_{\varepsilon \downarrow 0}h_\varepsilon
= (\tilde h)_0
= (\tilde h|_\gotA)_0
= (h|_\gotA)_0
= h_0
\]
where we used that $\gotA$ is dense in $D(\tilde h)$ and in $D(h)$.\hfill$\Box$

\ruimte

Next  we derive a crude  $L_2$ off-diagonal bound.

If $\ce$ is a Dirichlet form over $X$ then we denote by $S = S^{(ce)}$ the semigroup generated by 
the self-adjoint operator associated with $\ce$.

\begin{lemma} \label{llbs232}
Let $h\in \cc_l$ and $\lambda > 0$ with $h\leq\lambda\,l|_\gotA$.
If $A,B \subset X$ are measurable then 
\[
|(\psi,S^{(h_0)}_t \varphi)|
\leq e^{-(4 \lambda t)^{-1} d^{(l)}(A;B)^2} \, \|\psi\|_2 \, \|\varphi\|_2
\]
for all $t > 0$, $\varphi \in L_2(A)$ and $\psi \in L_2(B)$.
\end{lemma}
\proof\
For all $\varepsilon > 0$
the form $h_\varepsilon$ is a diffusion with 
$h_\varepsilon \leq (\lambda + \varepsilon) \,l$
and $D(h_\varepsilon) = D(l)$ by Corollary~\ref{clbs231.5}.
So
\[
\ci^{(h_\varepsilon)}_\psi(\varphi)
\leq \ci^{((\lambda + \varepsilon) l)}_\psi(\varphi)
= \ci^{(l)}_{(\lambda + \varepsilon)^{1/2} \psi}(\varphi)
\]
for all $\varphi,\psi \in D(l)\cap L_\infty(X)$ with $\varphi \geq 0$
by \cite{ERSZ2}, Proposition~3.2.
Therefore  
$|||\widehat \ci^{(h_\varepsilon)}_\psi||| \leq |||\widehat \ci^{(l)}_{(\lambda + \varepsilon)^{1/2} \psi}|||$
for all $\psi \in D(l)_{\rm loc} \cap L_\infty(X)= D(h_\varepsilon)_{\rm loc} \cap L_\infty(X)$.
So
 $(\lambda + \varepsilon)^{-1/2} d^{(l)}(A\,;B) \leq d^{(h_\varepsilon)}(A\,;B)$
for all measurable $A,B \subset X$.
 
The form $h_\varepsilon$ is regular by Corollary~\ref{clbs231.5}.
Then it follows from  Theorem~\ref{tdfd3.1}, that 
\[
|(\psi,S^{(h_\varepsilon)}_t \varphi)|
\leq e^{-(4 t)^{-1} d^{(h_\varepsilon)}(A;B)^2} \, \|\psi\|_2 \, \|\varphi\|_2
\leq e^{-(4 (\lambda + \varepsilon) t)^{-1} d^{(l)}(A;B)^2} \, \|\psi\|_2 \, \|\varphi\|_2
\]
for all $t > 0$, measurable $A,B \subset X$, $\varphi \in L_2(A)$ and $\psi \in L_2(B)$.
Since $\lim_{\varepsilon \downarrow 0} S^{(h_\varepsilon)}_t = S^{(h_0)}_t$
strongly, as a consequence of strong resolvent convergence, the lemma follows.\hfill$\Box$

\ruimte

To exploit the $L_2$ off-diagonal bound we need an estimate for $d^{(l)}(A\,;B)$.
One can obtain a simple estimate  if $l$ satisfies Condition~L.
{\bf Therefore from now on we always assume that $l$ satisfies {\rm Condition~L}.}

Recall that if  $A,B \subset X$ with $A \neq \emptyset$ and  $B \neq \emptyset$
then
$d^{(e)}(A\,;B) = \inf_{x\in A}\inf_{y\in B} \,d^{(e)}(x\,;y)$.

\begin{lemma} \label{las352}
\mbox{}
\begin{tabel}
\item \label{las352-1}
If $U,V \subset X$ are open and non-empty then 
$d^{(l)}(U\,;V) \geq d^{(e)}(U\,;V)$.
\item \label{las352-2}
If $K_1,K_2 \subset X$ are compact with $K_1 \cap K_2 = \emptyset$ then 
$d^{(l)}(K_1\,;K_2) > 0$.
\end{tabel}
\end{lemma}
\proof\
The first statement is easy.
In order to prove the second one, we use regularity of the metric space $X$.
There are non-empty open subsets $U_1,U_2$ and closed subsets $F_1,F_2$ such that 
$K_1 \subset U_1 \subset F_1$, $K_2 \subset U_2 \subset F_2$ and $F_1 \cap F_2 = \emptyset$.
Then 
\[
d^{(l)}(K_1\,;K_2)
\geq d^{(l)}(U_1\,;U_2)
\geq d^{(e)}(U_1\,;U_2)
\geq d^{(e)}(F_1\,;F_2) > 0
\;\;\; ,  \]
from which the lemma follows.\hfill$\Box$

\ruimte

We next  examine local approximations for the Markovian forms $h\in\cc_l$ defined by truncation.
We adopt the definition used for Dirichlet forms and demonstrate that the domination
property allows one to deduce many of the standard properties for the relaxations.

If $h \in \cc_l$ and $\Phi \in \gotA$ with $\Phi \geq 0$ then the quadratic form $h_\Phi$
is defined by $D(h_\Phi)=\gotA$ and 
\[
h_\Phi(\varphi) 
= h(\Phi\varphi,\varphi) - 2^{-1} h(\Phi,\varphi^2)
\]
for all $\varphi\in \gotA$.
This is well-defined since $\gotA$ is an algebra.
Moreover, $h_\Phi(\varphi) = \ci^{(h)}_\varphi(\Phi)$ if $h$ is a Dirichlet form.
Now the truncation $h_{\varepsilon\,\Phi}$ of the Dirichlet form $h_\varepsilon$
satisfies 
\[
0
\leq h_{\varepsilon\,\Phi}(\varphi)
= \ci^{(h_\varepsilon)}_\varphi(\Phi)
\leq \|\Phi\|_\infty \, h_\varepsilon(\varphi)
= \|\Phi\|_\infty \Big( h(\varphi) + \varepsilon \, l(\varphi) \Big)
\]
for all $\varphi\in \gotA$ and $\varepsilon>0$,
where we have used (\ref{easr1.11}).
But $h_{\varepsilon\,\Phi}(\varphi) = h_\Phi(\varphi)+\varepsilon\,l_\Phi(\varphi)$.
Hence in the limit $\varepsilon\downarrow 0$ one deduces that
\begin{equation}
0\leq h_\Phi(\varphi)\leq \|\Phi\|_\infty\,h(\varphi)\leq \lambda\,\|\Phi\|_\infty\,l(\varphi)
\label{eSas3;8}
\end{equation}
for all $\varphi\in \gotA$.
One deduces in a similar fashion that $h_\Phi$ satisfies the basic Markov property:
$\varphi\in \gotA$ implies 
$0 \vee \varphi\wedge \one\in \gotA$ and $h_\Phi(0 \vee \varphi\wedge \one) \leq h_\Phi(\varphi)$.
Moreover, $h_\Phi$ is strongly local.
Hence $h_\Phi \in \cc_l$.

Note that in general $h_\Phi$ 
is not a Dirichlet form as it is not necessarily closed and possibly not closable, even
if $h$ is a Dirichlet form.
The relaxation, or viscosity closure, $h_{\Phi\,0}=(h_\Phi)_0$ is, however, a regular Dirichlet form
by Lemma~\ref{llbs231} and the next lemma establishes that it is even a diffusion if $\supp \Phi$ is compact.

\begin{lemma} \label{llbs233}
Let  $h\in \cc_l$ and  $\Phi \in \gotA$ with $0\leq \Phi \leq 1$
and $\supp \Phi$ compact.
Then $h_{ \Phi \, 0}\in\cc_l$ and $h_{ \Phi \, 0}$ is a diffusion.

Moreover, $S^{(h_{\Phi \, 0})}$ leaves $L_2(\supp \Phi)$ invariant
and $S^{(h_{\Phi \, 0})}_t \varphi = \varphi$ for all $\varphi \in L_2((\supp \Phi)^{\rm c})$
and $t > 0$.
\end{lemma}
\proof\
For brevity write $S = S^{(h_{ \Phi \, 0})}$.

It follows from the foregoing that $h_\Phi \in\cc_l$.
Therefore  $h_{\Phi\,0}$ 
is a regular Dirichlet form by the first statement of Lemma~\ref{llbs231}
and $D(l)$ is dense in $ D(h_{\Phi\,0})$ by the second statement of the lemma.
Next we establish the localization properties of  $h_{ \Phi \, 0}$ and the semigroup $S$.

Since $l$ is regular there are $\chi_1,\chi_2,\ldots \in D(l) \cap C_c(X)$ such that 
$\lim_{n \to \infty} \chi_n = \one_{\supp \Phi}$ in $L_2(X)$ and for all 
$n \in \Ni$ there is a neighbourhood $U$ of $\supp \Phi$ with $\chi_n|_U = 1$.
Then 
\[
h_{ \Phi}(\chi_n) = h(\Phi,\chi_n)-2^{-1}h(\Phi,\chi_n^2)=0
\]
because $h$ is strongly local.
Moreover
\[
h_{\Phi \, 0}(\one_{\supp \Phi}) 
      \leq \liminf_{n \to \infty} h_{ \Phi}(\chi_n) = 0
\;\;\;.
\]
Therefore $\one_{\supp \Phi} \in D(h_{\Phi \, 0})$ and $h_{\Phi \, 0}(\one_{\supp \Phi}) = 0$.
Hence $S_t \one_{\supp \Phi} = \one_{\supp \Phi}$
for all $t > 0$ and $S$ leaves $L_2(\supp \Phi)$ invariant.
(Cf.\ \cite{ERSZ1} Lemma~6.1.)
Similarly, if $\varphi \in \gotA$ then 
$h_{\varepsilon \, \Phi}((\one - \chi_n) \varphi) =-2^{-1}h_\varepsilon (\Phi,(1-\chi_n)^2 \, \varphi^2)=0$ 
for all $n \in \Ni$ and $\varepsilon > 0$ by [BH]-locality of $h_\varepsilon$.
So $h_\Phi((\one - \chi_n) \varphi) = \lim_{\varepsilon \downarrow 0} h_{\varepsilon \, \Phi}((\one - \chi_n) \varphi) = 0$
and 
\[
h_{\Phi \, 0}(\one_{(\supp \Phi)^{\rm c}} \varphi)
      \leq \liminf_{n \to \infty} h_{\Phi}((\one - \chi_n) \varphi) = 0
\;\;\;.
\]
Now let $\varphi \in L_2((\supp \Phi)^{\rm c})$.
Since $\gotA$ is dense in $L_2(X)$ there are 
$\varphi_1,\varphi_2,\ldots \in \gotA$ such that 
$\lim_{n \to \infty} \varphi_n = \varphi$ in $L_2(X)$.
Then $\lim_{n \to \infty} \one_{(\supp \Phi)^{\rm c}} \, \varphi_n = \varphi$ in $L_2(X)$.
Since $h_{ \Phi \, 0}$ is lower semicontinuous one deduces that 
$h_{\Phi \, 0}(\varphi) 
       \leq \liminf_{n \to \infty} h_{\Phi \, 0}(\one_{(\supp \Phi)^{\rm c}} \varphi) = 0$.
So $\varphi \in D(h_{\Phi \, 0})$
and $h_{\Phi \, 0}(\varphi) = 0$.
Hence $S^{(h_{\Phi \, 0})}_t \varphi = \varphi$ for all $\varphi \in L_2((\supp \Phi)^{\rm c})$
and $t > 0$.

Finally we prove that $h_{ \Phi \, 0}$ is strongly local.
First, if ${h_\Phi} \leq\lambda\,l|_\gotA$ with $\lambda > 0$ 
then
\[
|(\psi,S_t \varphi)|
\leq e^{-(4 \mu t)^{-1} d^{(l)}(A\,;B)^2} \, \|\psi\|_2 \, \|\varphi\|_2
\]
for all measurable $A,B \subset X$, $t > 0$, $\varphi \in L_2(A)$ and $\psi \in L_2(B)$
by Lemma~\ref{llbs232} and (\ref{eSas3;8}), where $\mu = \lambda \, \|\Phi\|_\infty + 1$.
Secondly, let $\varphi,\psi \in D(h_{\Phi \, 0})$ with $\supp \varphi$ and $\supp \psi$ compact 
and suppose there exists a neighbourhood $U$ of $\supp \varphi$ such that 
$\psi|_U = 1$.
There exists a $\chi \in D(l) \cap C_c(X)$ such that $K = \supp \chi$ is compact
and $\chi|_{\supp \psi \cup \supp \Phi} = 1$.
Then $h_{ \Phi \, 0}(\chi) = 0$
since $\chi = \one_{\supp \Phi} + \chi \, \one_{(\supp \Phi)^{\rm c}}$.
Therefore $S_t \chi = \chi$ for all $t > 0$.
Then $(\psi,\varphi) = (\chi, \varphi) = (\chi, S_t \varphi)$ for all $t > 0$.
Hence 
\[
t^{-1} \, |(\psi, (I - S_t) \varphi)|
= t^{-1} |(\chi - \psi, S_t \varphi)|  
\leq  t^{-1} e^{- (4 \mu t)^{-1} d^{(l)}(K \backslash U ; \supp \varphi)^2} \, 
           \|\chi - \psi\|_2 \, \|\varphi\|_2
\]
for all $t > 0$.
But $d^{(l)}(K \backslash U ; \supp \varphi) > 0$ by Condition~L and Lemma~\ref{las352}.
So 
\[
|h_{ \Phi \, 0}(\psi,\varphi)|
= \lim_{t \downarrow 0} t^{-1} \, |(\psi, (I - S_t) \varphi)|
= 0
\;\;\; .  \]
Thus $h_{ \Phi \, 0}$ is strongly local. 
This completes the verification that $h_{ \Phi \, 0}$
is a diffusion and  $h_{ \Phi \, 0}\in\cc_l$.\hfill$\Box$

\ruimte

If $h \in \cc_l$ and $\Phi \in \gotA$ with $\supp \Phi$ compact and $0 \leq \Phi \leq 1$
write $\ci^{(\Phi)} = \ci^{(h_{ \Phi \, 0})}$ and $d^{(\Phi)} = d^{(h_{\Phi \, 0})}$.
Further let $H_{ \Phi \, 0}$ and $S^{(\Phi)}$ denote the operator and semigroup 
associated with the form $h_{ \Phi \, 0}$.
Moreover, let $H_0$ denote the operator associated with the form $h_0$.
It will be clear from the context which $h$ is involved.
Proposition~\ref{llbs233} establishes that $h_{\Phi\,0}$ is a localization of $h$: the corresponding
semigroup $S^\Phi$ leaves  $L_2(\supp\Phi)$ invariant.
Next we consider the distance corresponding to the generator of the restriction of 
$S^\Phi$ to $L_2(\supp\Phi)$ 
and its relation to the distance corresponding to $h_{\Phi\,0}$.

\begin{lemma} \label{llbs105}
Let $h\in\cc_l$ and $\Phi \in \gotA$ with $0 \leq \Phi \leq 1$ and $\supp \Phi$ compact.
Define the form $\hat h$ on $L_2(\supp\Phi)$ by 
\[
D(\hat h) = \{ \varphi|_{\supp \Phi} : \varphi \in D(h_{ \Phi \, 0}) \} 
\]
and 
\[
\hat h(\varphi|_{\supp \Phi}) = h_{ \Phi \, 0}(\varphi \one_{\supp \Phi})
\]
for all $\varphi \in D(h_{ \Phi \, 0})$.
Then 
\[
d^{(\Phi)}(A\,;B) 
= d^{(\hat h)}(A \cap \supp \Phi \,;B \cap \supp \Phi) 
\]
for all measurable $A,B \subset X$.
\end{lemma}
\proof\
Note that the form $\hat h$ is a conservative [BH]-local Dirichlet form.
Moreover, $D(\hat h)_{\rm loc} = D(\hat h)$.
If $\varphi,\psi \in D(h_{ \Phi \, 0}) \cap L_\infty$ then
\[
\ci^{(\Phi)}_{\psi \one_{\supp \Phi}}(\varphi)
= h_{ \Phi \, 0}(\psi \one_{\supp \Phi} \, \varphi, \psi \one_{\supp \Phi}) 
     - 2^{-1} h_{ \Phi \, 0}(\varphi, \psi^2 \one_{\supp \Phi})
= \ci^{(\hat h)}_{\psi|_{\supp \Phi}}(\varphi|_{\supp \Phi})
\]
since $h_{ \Phi \, 0}(\one_{(\supp \Phi)^{\rm c}} \varphi) = 0$.
So $|||\widehat \ci^{(\Phi)}_{\psi \one_{\supp \Phi}}||| = |||\widehat \ci^{(\hat h)}_{\psi|_{\supp \Phi}}|||$
for all $\psi \in D(h_{\Phi \, 0})_{\rm loc} \cap L_\infty$.

Let $M \in [0,d^{(\Phi)}(A\,;B)] \cap \Ri$ and $\varepsilon > 0$.
By \cite{ERSZ2}, Lemma~4.2.III, there exists a $\psi \in D_0(h_{ \Phi \, 0})$ such that 
$\psi|_B = 0$ and $\psi|_A \geq M - \varepsilon$.
Then $\psi|_{\supp \Phi} \in D_0(\hat h)$.
Therefore 
\[
M - \varepsilon 
\leq d_\psi(A\,;B) 
\leq d_{\psi|_{\supp \Phi}}(A \cap \supp \Phi\,;B \cap \supp \Phi)
\leq d^{(\hat h)}(A \cap \supp \Phi \,; B \cap \supp \Phi)
\;\;\; .  \]
So $d^{(\Phi)}(A\,;B) \leq d^{(\hat h)}(A \cap \supp \Phi \,; B \cap \supp \Phi)$.

Conversely, let $M \in [0,d^{(\hat h)}(A \cap \supp \Phi \,; B \cap \supp \Phi)] \cap \Ri$ and $\varepsilon > 0$.
By another application of Lemma~4.2.III in \cite{ERSZ2} there exists a $\tilde \psi \in D_0(\hat h)$ such that 
$\tilde \psi|_{B \cap \supp \Phi} = 0$ and $\tilde \psi|_{A \cap \supp \Phi} \geq M - \varepsilon$.
There exists a $\psi \in D(h_{\Phi \, 0})$ such that $\psi|_{\supp \Phi} = \tilde \psi$.
Then $\psi \one_{\supp \Phi} \in D_0(h_{\Phi \, 0})$.
Moreover, $\widehat \ci^{(\Phi)}_{\one_{(\supp \Phi)^{\rm c} \cap A}} = 0$ by the last part of Lemma~\ref{llbs233}.
Let $\tau = (M - \varepsilon) \one_{(\supp \Phi)^{\rm c} \cap A} + \psi \one_{\supp \Phi}$.
Then $\tau \in D_0(h_{\Phi \, 0})_{\loc}$ by \cite{ERSZ2}, Lemma~3.3.I, 
 $\tau|_B = 0$ and $\tau|_A \geq M - \varepsilon$.
Therefore 
\[
M-\varepsilon 
\leq d_\tau(A\,;B)
\leq d^{(\Phi)}(A\,;B)
\]
and $d^{(\hat h)}(A \cap \supp \Phi \,; B \cap \supp \Phi) \leq d^{(\Phi)}(A\,;B)$.\hfill$\Box$

\ruimte

Now we can make the first key deduction in the proof of the theorem.
One can apply the results of Hino--Ram\'{\i}rez \cite{HiR} to the form $\hat h$ in Lemma~\ref{llbs105}.

\begin{cor} \label{clbs105.5}
Let $h\in\cc_l$ and $\Phi \in \gotA$ with $0 \leq \Phi \leq 1$ and $\supp \Phi$ compact.
Let $A,B \subset X$ be  measurable with $A,B \subseteq \supp \Phi$.
Then 
\[
\lim_{t \downarrow  0} t \log (\one_A, S^{(\Phi)}_t \one_B)
= - 4^{-1} d^{(\Phi)}(A\,;B)^2
\;\;\; .  
\]
\end{cor}
\proof\
This follows from Lemma~\ref{llbs105} and \cite{HiR}, Theorem~1.1.\hfill$\Box$

\ruimte

The next idea is to take a limit over the localizations $\Phi$ as $\Phi$ increases
monotonically to the identity function.
This involves analyzing both the limit of the semigroups $S^{(\Phi)}$
and of the distances $d^{(\Phi)}$.
It is not difficult using arguments of monotonicity to deduce that the strong limit
of the semigroups $S^{(\Phi)}$ and the pointwise limit of the distances $d^{(\Phi)}$
exist. 
But we also have to identify the limits with the semigroup and distance corresponding to the relaxation
and to  control simultaneously the small time behaviour.
The key to this analysis is the observation that the associated wave equations have a
finite speed of propagation (see \cite{ERSZ1}, Proposition~3.2).

\smallskip

If $A \subset X$ with $A \neq \emptyset$ and $r > 0$ set 
\[
\widetilde B^{(e)}(A\,;r) 
= \{ \,x \in X : d^{(e)}(x\,;A) \leq r \,\} 
\;\;\; ,  \]
where $d^{(e)}(x\,;A) = \inf \{ |x-a| : a \in A \} $.
Note that $\widetilde B^{(e)}(A\,;r)$ is closed since $x \mapsto d^{(e)}(x\,;A)$ is continuous.
Since $h_\Phi \leq \|\Phi\|_\infty \, h|_\gotA$ by (\ref{eSas3;8})
one has $h_{\Phi \, 0} \leq \|\Phi\|_\infty \, (h|_\gotA)_0 = \|\Phi\|_\infty \, h_0$, so
clearly $D(h_0) \subset D(h_{ \Phi \, 0})$ if $\Phi \in D(l) \cap L_\infty$ with $\Phi \geq 0$.

\begin{prop} \label{plbs106}
Let  $h\in \cc_l$, $\lambda > 0$  and  $\Phi \in \gotA$ and suppose that $h\leq\lambda\,l|_\gotA$, $0\leq \Phi \leq 1$
and $\supp \Phi$ is compact.
Further let $A \subset X$ measurable, $\Omega \subset X$ open with $\emptyset \neq \overline A \subset \Omega$
and suppose $\Phi|_\Omega = 1$.
Then 
\begin{equation}
\cos(t H_0^{1/2}) \varphi
= \cos(t H_{ \Phi \, 0}^{1/2}) \varphi
\label{eplbs106;2}
\end{equation}
for all $t \in \Ri$ with $|t| \leq \lambda^{-1/2} d^{(e)}(A\,;\Omega^{\rm c})$
and $\varphi \in L_2(A)$.
Moreover, if $\psi \in D(h_{ \Phi \, 0})$ and $\supp \psi \subset \Omega$ then 
$\psi \in D(h_0)$ and $h_0(\psi) = h_{ \Phi \, 0}(\psi)$.
\end{prop}

The proof of the proposition relies on two lemmas.

\begin{lemma} \label{llbs107}
Let $h\in\cc_l$ be a diffusion, $\lambda > 0$ with $h\leq\lambda\,l|_\gotA$, 
and let $H$ denote the corresponding positive self-adjoint operator.
If $A \subset X$ is open with $A \neq \emptyset$ then 
\[
\cos(t H^{1/2}) L_2(A)
\subset L_2(\widetilde B^{(e)}(A\,; \lambda^{1/2} \, |t|) )
\]
and 
\[
(t H^{1/2})^{-1} \sin(t H^{1/2}) L_2(A)
\subset L_2(\widetilde B^{(e)}(A\,; \lambda^{1/2} \, |t|) )
\]
for all $t \in \Ri \backslash \{ 0 \} $.
\end{lemma}

Here and in the sequel the operator formally denoted 
by $(t H^{1/2})^{-1} \sin(t H^{1/2})$ is properly defined by spectral theory, even if $H = 0$.

\ruimte

\noindent
\proof\
It follows from Lemmas~\ref{llbs232} and \ref{las352}.\ref{las352-1} that 
\[
|(\psi, S^{(h)}_t \varphi)|
\leq e^{- (4 \lambda t)^{-1} d^{(e)}(A\,;B)^2} \|\psi\|_2 \, \|\varphi\|_2
\]
for all $t > 0$, non-empty open $B \subset X$, $\varphi \in L_2(A)$ and $\psi \in L_2(B)$.
Therefore
\[
(\psi, \cos(t H^{1/2}) \varphi) = 0
\]
for all non-empty open $B \subset X$, $\varphi \in L_2(A)$, $\psi \in L_2(B)$
and $t \in \Ri$ with $|t| \leq \lambda^{-1/2} \, d^{(e)}(A\,;B)$ 
by \cite{ERSZ1}, Lemma~3.3.
Hence if $r > 0$, $\psi \in C_c(\widetilde B^{(e)}(A\,; r)^{\rm c})$ and $\varphi \in L_2(A)$
then 
$(\psi, \cos(t H^{1/2}) \varphi) = 0$ 
for all $|t| \leq \lambda^{-1/2} \, r$.
This implies the first identity.
The second identity  follows from the observation that
\[
(\psi, (t H^{1/2})^{-1} \sin(t H^{1/2}) \varphi)
= t^{-1} \int_0^t ds \, (\psi, \cos(s H^{1/2}) \varphi)
\]
combined with  the first identity.\hfill$\Box$

\begin{lemma} \label{llbs108}
Let $h\in\cc_l$ be a diffusion 
and let $H$ denote the corresponding positive self-adjoint operator.
If $A \subset X$ is open then 
$\overline{ D(H) \cap L_2(A) } =\makebox[0pt]{\raisebox{6mm}{}} L_2(A)$.
\end{lemma}
\proof\
We may assume that $A \neq \emptyset$.
Let $\lambda > 0$ and suppose that $h\leq\lambda\,l|_\gotA$.
Let $\varphi \in C_c(A)$ with $\varphi \neq 0$.
For all $t > 0$ set $\varphi_t = (t^2 H)^{-1} \sin^2(t H^{1/2}) \varphi$.
Then $\varphi_t \in D(H)$ and 
$\varphi_t \in L_2(\widetilde B^{(e)}(\supp \varphi\,; 2 \lambda^{1/2} \, |t|) )$
for all $t > 0$ by Lemma~\ref{llbs107}.
Hence $\varphi_t \in L_2(A)$ if $t > 0$ is small enough.
Finally, $\lim_{t \downarrow 0} \varphi_t = \varphi$ by spectral theory.
So $C_c(A) \subset \overline{ D(H) \cap L_2(A) }$ 
and the lemma follows.\hfill$\Box$

\ruimte

The conclusion of the last lemma is very useful
since it  shows that the operator  domain contains abundant functions with 
compact support.

\ruimte

\noindent
{\bf Proof of Proposition~\ref{plbs106}\hspace{5pt}}\
Fix $\varepsilon > 0$.
We begin by comparing the actions of the operators $H_{\varepsilon}$ and $H_{\Phi\,\varepsilon}$
corresponding to the diffusions $h_{\varepsilon}$ and
$h_{\Phi\,\varepsilon}$ associated with $h$ and $h_\Phi$
(see Corollary~\ref{clbs231.5}).

First assume that $A$ is an open subset of $X$ and let
 $r \in \langle0, (\lambda + \varepsilon)^{-1/2} \, d^{(e)}(A\,;\Omega^{\rm c}) \rangle$.
Fix  $\varphi\in L_2(A)$. 
Then it follows from Lemma~\ref{llbs107} applied to $h_\varepsilon$  that 
\[
\cos(t H_{ \varepsilon}^{1/2}) \varphi
\in L_2(\widetilde B^{(e)}(A\,; (\lambda + \varepsilon)^{1/2} \, |t|) )
\subset L_2(\widetilde B^{(e)}(A\,; (\lambda + \varepsilon)^{1/2} \, r) )
\subset L_2(\Omega)
\]
for all $t \in \langle-r,r\rangle$.

If in addition $\varphi \in D(h_{ \, \varepsilon})=D(l)$ then 
$\cos(t H_{ \varepsilon}^{1/2}) \varphi \in D(l)$.
Thus if one sets $\varphi_n=(-n) \vee \cos(t H_{ \varepsilon}^{1/2}) \varphi \wedge n$ for all $n \in \Ni$ one has
$\varphi_n\in D(l)\cap L_\infty(X) = \gotA$
and $\supp \varphi_n \subset \widetilde B^{(e)}(A\,; (\lambda + \varepsilon)^{1/2} \, r)$.
Since $l$ is regular there exists a $\chi \in \gotA$ with 
$0 \leq \chi \leq 1$, $\supp \chi \subset \Omega$ and $\chi|_{\widetilde B^{(e)}(A; (\lambda + \varepsilon)^{1/2} \, r)} = 1$.
Let $\psi \in \gotA$.
Then $\chi\,\psi\in \gotA$ and $h_{\Phi \, \varepsilon}(\psi, \varphi_n) = h_{\Phi \, \varepsilon}(\chi \, \psi,  \varphi_n) $
by locality.
But
\begin{eqnarray*}
h_{\Phi \, \varepsilon}( \chi \,\psi,  \varphi_n) 
&=&h_\Phi(\chi \, \psi,  \varphi_n) + \varepsilon\,l(\chi \, \psi,  \varphi_n) \\[5pt]
&=&h(\chi \, \psi,  \varphi_n) + \varepsilon\,l(\chi \, \psi,  \varphi_n)
=h_{\varepsilon}(\chi \, \psi,  \varphi_n)
\end{eqnarray*}
by the definition of $h_\Phi$ and the assumption $\Phi|_\Omega = 1$ .
Therefore
\[
h_{\Phi \, \varepsilon}(\psi, \varphi_n)
=h_{\varepsilon}(\chi \, \psi,  \varphi_n)
=h_{\varepsilon}( \psi,  \varphi_n)
\]
where the second identity follows by another use of locality.
Then  the limit $n\to\infty$ gives
\[
h_{\Phi\,\varepsilon}( \psi,  \cos(t H_{ \varepsilon}^{1/2}) \varphi)
=h_\varepsilon(\psi, \cos(t H_{ \varepsilon}^{1/2}) \varphi)
\]
for all $\psi\in \gotA$, $\varphi\in D(h_\varepsilon)\cap L_2(A)$
and $t \in \langle-r,r\rangle$ by \cite{FOT}, Theorem~1.4.2(iii).

Let $\varphi \in D(H_{ \varepsilon}) \cap L_2(A)$ and $t \in \langle-r,r\rangle$.
Because $\gotA$ is a core for $h_{\Phi \, \varepsilon}$ one has 
\[
h_{ \Phi \, \varepsilon}(\psi, \cos(t H_{ \varepsilon}^{1/2}) \varphi)
= h_{ \varepsilon}(\psi, \cos(t H_{ \varepsilon}^{1/2}) \varphi)
= (\psi, \cos(t H_{ \varepsilon}^{1/2}) H_{ \varepsilon} \varphi)
\]
for all $\psi \in D(h_{ \Phi \, \varepsilon})$.
Since the form $h_{ \Phi \, \varepsilon}$ is closed it follows that
 $\cos(t H_{ \varepsilon}^{1/2}) \varphi \in D(H_{ \Phi \, \varepsilon})$ and 
\begin{equation}
H_{ \Phi \, \varepsilon} \cos(t H_{ \varepsilon}^{1/2}) \varphi = \cos(t H_{ \varepsilon}^{1/2}) H_{ \varepsilon} \varphi
\label{eplbs106;1}
\end{equation}
for all $t \in \langle-r,r\rangle$ and $\varphi \in D(H_{ \varepsilon}) \cap L_2(A)$.

Next let $\varphi \in D(H_{ \varepsilon}) \cap L_2(A)$ and 
for all $t \in \Ri$ define 
\[
\chi_t = \cos(t H_{ \varepsilon}^{1/2}) \varphi - \cos(t H_{ \Phi \, \varepsilon}^{1/2}) \varphi
\;\;\; .  \]
Clearly $\|\chi_t\|_2 \leq 2\,\|\varphi\|_2$ for all $t \in \Ri$.
Let $\psi \in D(H_{ \Phi \, \varepsilon})$. 
Then $\tau \colon t \mapsto (\psi,\chi_t)$ is twice differentiable and $\tau(0) = \tau'(0) = 0$.
Moreover,
\begin{eqnarray*}
(\psi, \chi_t)
& = & \int_0^t dt_1 \int_0^{t_1} dt_2 \, \tau''(t_2)  \\[5pt]
& = & \int_0^t dt_1 \int_0^{t_1} dt_2 \, \Big((H_{ \Phi \, \varepsilon} \psi, \cos(t_2 H_{ \Phi \, \varepsilon}^{1/2}) \varphi)
    - (\psi, \cos(t_2 H_{ \varepsilon}^{1/2}) H_{ \varepsilon} \varphi)\Big)  \\[5pt]
& = & - \int_0^t dt_1 \int_0^{t_1} dt_2 \, (H_{ \Phi \, \varepsilon} \psi, \chi_{t_2})
\end{eqnarray*}
for all $t \in \langle-r,r\rangle$, where we used (\ref{eplbs106;1}) in the last step.
Now assume that $\psi$ is a bounded vector for $H_{ \Phi \, \varepsilon}$, i.e., 
$\psi \in \bigcap_{n=1}^\infty D(H_{ \Phi \, \varepsilon}^n)$ and 
there exists a $b > 0$ such that $\|H_{ \Phi \, \varepsilon}^n \psi\|_2 \leq b^n$ for all $n \in \Ni$.
Then by iteration one deduces that 
\begin{eqnarray*}
|(\psi, \chi_t)|
& = & \Big| \int_0^t dt_1 \ldots \int_0^{t_{2n-1}} dt_{2n} \, (H_{ \Phi \, \varepsilon}^n \psi, \chi_{t_{2n}}) \Big|  \\[5pt]
& \leq & 2 \,|t|^{2n} \, (2n)!^{-1} \|H_{ \Phi \, \varepsilon}^n \psi\|_2 \, \|\varphi\|_2  
\leq 2\, |t|^{2n} \, (2n)!^{-1} b^n \, \|\varphi\|_2  
\end{eqnarray*}
for all $n \in \Ni$ and $t \in \langle-r,r\rangle$.
Taking the limit $n \to \infty$ it follows that $(\psi,\chi_t) = 0$ for all 
bounded vectors $\psi$ for $H_{ \Phi \, \varepsilon}$ and $t \in \langle-r,r\rangle$.
Since the bounded vectors for $H_{ \Phi \, \varepsilon}$ are dense in $L_2$ by spectral theory one deduces that 
$\cos(t H_{ \Phi \, \varepsilon}^{1/2}) \varphi = \cos(t H_{ \varepsilon}^{1/2}) \varphi$
for all $t \in \langle-r,r \rangle$ and $\varphi \in D(H_{ \varepsilon}) \cap L_2(A)$.
Hence 
\[
\cos(t H_{ \Phi \, \varepsilon}^{1/2}) \varphi = \cos(t H_{ \varepsilon}^{1/2}) \varphi
\]
for all $t \in \langle-r,r \rangle$ and $\varphi \in L_2(A)$
by Lemma~\ref{llbs108}, since $A$ is open.

Let $\varphi \in L_2(A)$.
Then $\cos(t H_{ \varepsilon}^{1/2}) \varphi = \cos(t H_{ \Phi \, \varepsilon}^{1/2}) \varphi$
for all $\varepsilon > 0$ and for all $t$ with 
$|t| < (\lambda + \varepsilon)^{-1/2} \, d^{(e)}(A\,;\Omega^{\rm c})$.
Therefore
\[
\cos(t H_0^{1/2}) \varphi
= \lim_{\varepsilon \downarrow 0} \cos(t H_{ \varepsilon}^{1/2}) \varphi
= \lim_{\varepsilon \downarrow 0} \cos(t H_{ \Phi \, \varepsilon}^{1/2}) \varphi  
= \cos(t H_{ \Phi \, 0}^{1/2}) \varphi
\]
for all $|t| < \lambda^{-1/2} \, d^{(e)}(A\,;\Omega^{\rm c})$
by strong resolvent convergence of $H_{ \varepsilon}$ to $H_0$ and of 
$H_{ \Phi \, \varepsilon}$ to $H_{ \Phi \, 0}$
(see \cite{RS1}, Theorem~VIII.20 and Corollary~\ref{clbs231.5}).
Then the first statement of the proposition follows if $A$ is open.

If $A$ is not open then for all $\delta \in \langle0,d^{(e)}(A\,;\Omega^{\rm c}) \rangle$
one can apply the above to the open set
$B^{(e)}(A\,;\delta) = \{ x \in X : d^{(e)}(x\,;A) < \delta \} $ and deduce that 
(\ref{eplbs106;2}) is valid for all $\varphi \in L_2(A)$ and 
$|t| < \lambda^{-1/2} \, d^{(e)}(B^{(e)}(A\,;\delta) \,;\Omega^{\rm c})$.
Since $\lim_{\delta \downarrow 0} d^{(e)}(B^{(e)}(A\,;\delta) \,;\Omega^{\rm c}) = d^{(e)}(A\,;\Omega^{\rm c})$
by the triangle inequality
the first statement of the proposition follows for general $A$.

The last statement of the proposition follows since
\[
D(h_0)
= \{ \,\varphi \in L_2(X) : \sup_{t \in \langle0,1]} 2 \,t^{-2}(\varphi, (I - \cos(t H_0^{1/2})) \varphi) < \infty\, \}
\]
and $h_0(\varphi) = \lim_{t \downarrow 0} 2 \,t^{-2}(\varphi, (I - \cos(t H_0^{1/2})) \varphi)$
for all $\varphi \in D(h_0)$, with a similar expression for $h_{ \Phi \, 0}$.\hfill$\Box$

\ruimte

The proposition has three very useful corollaries.
The first corollary establishes the first statement of Theorem~\ref{tlbs121}.

\begin{cor} \label{clbs106.4}
If $h \in \cc_l$ then $h_0$ is a diffusion.
In particular, $h_0 \in \cc_l$.
\end{cor}
\proof\
By Lemma~\ref{llbs231} it remains to show that $h_0$ is strongly local.
Let $\varphi,\psi \in D(h_0)$ with $\supp \varphi$ and $\supp \psi$ compact and $\psi$ is
constant on a neighbourhood of $\supp \varphi$.
Since $l$ is regular there exist $\Phi \in \gotA$ and an open $\Omega \subset X$ such that 
$\supp (\varphi \pm \psi) \subset \Omega$, $0 \leq \Phi \leq 1$, $\supp \Phi$ is compact and 
$\Phi|_\Omega = 1$.
Then $h_{\Phi \, 0}(\psi,\varphi) = 0$ by Lemma~\ref{llbs233}.
But $h_0(\varphi \pm \psi) = h_{\Phi \, 0}(\varphi \pm \psi)$ by the last part of
Proposition~\ref{plbs106}.
Therefore $h_0(\psi,\varphi) = h_{\Phi \, 0}(\psi,\varphi) = 0$ by polarization
and symmetry.\hfill$\Box$

\ruimte

\begin{cor} \label{clbs106.5}
Adopt the assumptions of Proposition~$\ref{plbs106}$.
If  $\psi \in D(h_{ \Phi \, 0}) \cap L_\infty$ with $\supp \psi \subset \Omega$
then $|||\widehat \ci^{(h_0)}_\psi||| = |||\widehat \ci^{(\Phi)}_\psi|||$.
\end{cor}
\proof\
Since $l$ is regular there exists a
$\tau \in D(l) \cap L_\infty$ such that $0 \leq \tau \leq 1$, $\tau|_{\supp \psi} = 1$ and 
$\supp \tau \subset \Omega$.
Let $\varphi \in D(h_{ \Phi \, 0})_{\rm loc} \cap L_{\infty,c}$.
Then $\varphi \, \tau \in D(h_{ \Phi \, 0}) \cap L_\infty$ and $\varphi \, \tau \in D(h_0)$ 
by the last part of Proposition~\ref{plbs106}.
Hence by strong locality and Proposition~\ref{plbs106} one has
\begin{eqnarray*}
\widehat \ci^{(\Phi)}_\psi(\varphi)
= \ci^{(\Phi)}_\psi(\varphi \, \tau)
& = & h_{ \Phi \, 0}(\psi \, \varphi \, \tau, \psi) - 2^{-1} h_{ \Phi \, 0}(\varphi \, \tau, \psi^2)  \\[5pt]
& = & h_0(\psi \, \varphi \, \tau, \psi) - 2^{-1} h_0(\varphi \, \tau, \psi^2)  
= \widehat \ci^{(h_0)}_\psi(\varphi \, \tau)
\;\;\; .
\end{eqnarray*}
So 
\[
|\widehat \ci^{(\Phi)}_\psi(\varphi)|
\leq |||\widehat \ci^{(h_0)}_\psi||| \, \|\varphi \, \tau\|_1
\leq |||\widehat \ci^{(h_0)}_\psi||| \, \|\varphi\|_1
\]
and $|||\widehat \ci^{(\Phi)}_\psi||| \leq |||\widehat \ci^{(h_0)}_\psi|||$.
The opposite inequality is similar.\hfill$\Box$

\ruimte

The third corollary establishes that the semigroups $S^{(\Phi)}$ converge strongly
to the semigroup associated with the relaxation $h_0$.

\begin{cor} \label{clbs109}
Adopt the assumptions of Proposition~$\ref{plbs106}$.
Then
\[
\|S^{(h_0)}_t \varphi - S^{(\Phi)}_t \varphi\|_2
\leq 2 \,e^{-(4t)^{-1} r^2} \, \|\varphi\|_2
\]
for all $t > 0$ and $\varphi \in L_2(A)$, where
$r = \lambda^{-1/2} \, d^{(e)}(A\,; \Omega^{\rm c})$.
\end{cor}
\proof\
This follows from the identity
\begin{equation}
S^{(h_0)}_t \varphi = (\pi t)^{-1/2} \int_0^\infty ds\, e^{-(4t)^{-1} s^2} \cos(s H_0^{1/2}) \varphi
\label{eas3.297}
\end{equation}
and Proposition~\ref{plbs106}.\hfill$\Box$

\ruimte

Corollary~\ref{clbs109} gives good control over the semigroups $S^{(\Phi)}$ as
$\Phi\to 1$ and indirectly it gives control over the  distances.
We first deduce that the distances $d^{(\Phi)}(A\,;B)$ become independent of the choice of $\Phi$
as $\Phi\to 1$ with $A$ and $B$ fixed.

\begin{lemma} \label{llbs111}
Let  $h\in \cc_l$, $\lambda > 0$ 
and  $\Phi,\Psi \in \gotA$ 
and assume that  with $h\leq\lambda\,l|_\gotA$, $0 \leq \Phi,\Psi \leq 1$ and $\supp \Phi$, $\supp \Psi$ compact.
Further let $A,B \subset X$ measurable, $\Omega \subset X$ open with $\emptyset \neq \overline A \subset \Omega$,
$B \subset \supp \Phi \cap \supp \Psi$
and suppose $\Phi|_\Omega = \Psi|_\Omega = 1$.
If
\[
d^{(\Phi)}(A\,;B) \vee d^{(\Psi)}(A\,;B) \leq \lambda^{-1/2} d^{(e)}(A\,;\Omega^{\rm c})
  \]
then $d^{(\Phi)}(A\,;B) = d^{(\Psi)}(A\,;B)$.
\end{lemma}
\proof\
It follows from Proposition~\ref{plbs106} that 
\[
\cos(t H_{ \Phi \, 0}^{1/2}) \varphi
= \cos(t H_0^{1/2}) \varphi
= \cos(t H_{\Psi \, 0}^{1/2}) \varphi
\]
for all $t \in \Ri$ with $|t| \leq \lambda^{-1/2} d^{(e)}(A\,;\Omega^{\rm c})$
and $\varphi \in L_2(A)$.
Hence it follows as in the proof of Corollary~\ref{clbs109} that 
\[
\|S^{(\Phi)}_t \varphi - S^{(\Psi)}_t \varphi\|_2
\leq 2 \,e^{-(4t)^{-1} r^2} \, \|\varphi\|_2
\]
for all $t > 0$ and $\varphi \in L_2(A)$, where
$r = \lambda^{-1/2} \, d^{(e)}(A\,; \Omega^{\rm c})$.
In particular, 
\[
|(\one_B, S^{(\Phi)}_t \one_A) - (\one_B, S^{(\Psi)}_t \one_A)|
\leq 2 \,e^{-(4t)^{-1} r^2} \, |A|^{1/2} \, |B|^{1/2}
\]
and 
\begin{eqnarray*}
(\one_B, S^{(\Phi)}_t \one_A)
& \leq & (\one_B, S^{(\Psi)}_t \one_A) + 2 \,e^{-(4t)^{-1} r^2} \, |A|^{1/2} \, |B|^{1/2}  \\[5pt]
& \leq & e^{-(4t)^{-1} d^{(\Psi)}(A;B)^2} \, |A|^{1/2} \, |B|^{1/2} + 2\, e^{-(4t)^{-1} r^2} \, |A|^{1/2} \, |B|^{1/2}  \\[5pt]
& \leq & 3 \,|A|^{1/2} \, |B|^{1/2} \, e^{-(4t)^{-1} d^{(\Psi)}(A;B)^2}
\end{eqnarray*}
for all $t > 0$, where we used the Davies--Gaffney bounds of Theorem~\ref{tdfd3.1} in the second step.
It follows by Corollary~\ref{clbs105.5} and Lemma~\ref{llbs110} that 
$d^{(\Phi)}(A\,;B) \geq d^{(\Psi)}(A\,;B)$.
 The converse inequality is valid by a similar argument.
Hence the distances are equal.
\hfill$\Box$

\ruimte

Identification of the limit of the distances corresponding to the localizations
relies  on the  construction of  certain cut-off functions and this construction
is dependent on the topological assumption. 
The subsequent argument follows similar reasoning in \cite{BM}, Section~3, and \cite{Stu4}, Appendix~A.

\begin{lemma} \label{las345}

Let $y \in X$ and define $\psi \colon X \to [0,\infty\rangle$ by
$\psi(x) = d^{(e)}(x\,;y)$.
Then $\psi \wedge N \in D_0(l)$ for all $N \in \Ni$.
\end{lemma}
\proof\
Since $X$ is separable there are $y_1,y_2,\ldots \in X$ such that 
$ \{ y_k : k \in \Ni \} $ dense is in $X$.
Let $n \in \Ni$.
Then $X = \bigcup_{k=1}^\infty B^{(e)}(y_k \,; n^{-1})$.
Let $k \in \Ni$. By definition of $d^{(e)}$ there exists a 
$\psi_{nk} \in D_0(l) \cap C_{\rm b}(X)$ such that 
$\psi_{nk}(y) - \psi_{nk}(y_k) \geq d^{(e)}(y\,;y_k) - n^{-1}$.
Then 
\begin{eqnarray*}
\hspace{-5.5pt}\psi_{nk}(x) + d^{(e)}(x\,;y)
& \leq & \psi_{nk}(x) - \psi_{nk}(y_k) + \psi_{nk}(y_k) + d^{(e)}(y\,;y_k)
         - d^{(e)}(y\,;y_k) + d^{(e)}(x\,;y)  \\[5pt]
& \leq & d^{(e)}(x\,;y_k) + \psi_{nk}(y) + n^{-1} + d^{(e)}(x\,;y_k)  
 \leq  \psi_{nk}(y) + 3 \, n^{-1}
\end{eqnarray*}
for all $x \in B(y_k \,; n^{-1})$.
Set $\tilde \psi_{nk} = (\psi_{nk}(y) - \psi_{nk}) \vee 0$.
Then $\tilde \psi_{nk} \in D_0(l)$.
Since $\psi_{nk}(y) - \psi_{nk}(x) \leq d^{(e)}(x\,;y)$ for all $x \in X$ one 
has 
\[
0 \leq \tilde \psi_{nk}(x) \leq d^{(e)}(x\,;y)
\]
for all $x \in X$.
Moreover, if $x \in B(y_k \,; n^{-1})$ then 
\[
\tilde \psi_{nk}(x) 
\geq \psi_{nk}(y) - \psi_{nk}(x) 
\geq d^{(e)}(x\,;y) - 3 \, n^{-1}
\;\;\; .  \]
Hence $\psi = \sup_{n \in \Ni} \sup_{k \in \Ni} \tilde \psi_{nk}$
and the lemma follows from Lemma~\ref{las344}.\hfill$\Box$

\begin{cor} \label{cas346}
If $K \subset X$ is compact and $\varepsilon > 0$ then 
there exists a $\psi \in D(l) \cap C_c(X)$ such that $0 \leq \psi \leq 1$,
$\psi|_K = 1$ and $|||\widehat \ci^{(l)}_\psi||| \leq \varepsilon$.
\end{cor}
\proof\
We may assume that $K \neq \emptyset$.
Let $y \in K$. 
There exists an $M > 0$ such that $d^{(e)}(x\,;y) \leq M$ for all $x \in K$.
Define $\psi(x) = 0 \vee (N - \varepsilon^{1/2} \, d^{(e)}(x\,;y)) \wedge 1$
where $N = 1 + \varepsilon^{1/2} \, M$.
Then the corollary follows from Lemma~\ref{las345}.\hfill$\Box$

\ruimte

At this point we are prepared to prove the second statement in Theorem~\ref{tlbs121}.
The first step in the proof is to identify the limit of the distances corresponding to the localizations
as $\Phi$ increases to the identity.

\begin{lemma} \label{llbs112}
Let $h \in \cc_l$.
Fix $x \in X$.
Let $\Phi_1,\Phi_2,\ldots \in \gotA$ be such that 
$0 \leq \Phi_1 \leq \Phi_2 \leq \ldots \leq 1$, suppose $\Phi_n|_{B^{(e)}(x;n)} = 1$ and
$\supp \Phi_n$ is compact for all $n \in \Ni$.
Let $A,B \subset X$ be relatively compact  and measurable.
Then $d^{(h_0)}(A\,;B) = d^{(\Phi_n)}(A\,;B)$ for all large $n \in \Ni$.
\end{lemma}
\proof\
One has $d^{(h_0)}(A\,;B) \leq d^{(\Phi_{n+1})}(A\,;B) \leq d^{(\Phi_n)}(A\,;B)$ for all $n \in \Ni$ by 
(\ref{easr1.12}), Corollary~\ref{cas346} and 
Proposition~5.1 of \cite{ERSZ2}.
Hence 
\[
d^{(h_0)}(A\,;B) \leq \inf_{n \in \Ni} d^{(\Phi_n)}(A\,;B)
= \lim_{n \to \infty} d^{(\Phi_n)}(A\,;B)
\;\;\; .  
\]
Let $M \in [0,\infty\rangle$ and suppose that $M \leq \inf_{n \in \Ni} d^{(\Phi_n)}(A\,;B)$.
Let $\varepsilon \in \langle0,1]$.
By Corollary~\ref{cas346} there exists a $\chi \in D(l) \cap C_c(X)$ such that $0 \leq \chi \leq 1$,
$\chi|_{A \cup B} = 1$ and $|||\widehat \ci^{(l)}_\chi||| \leq \varepsilon^2$.
There is an $n \in \Ni$ such that $\supp \chi \subset B^{(e)}(x\,;n)$.
Then $M \leq d^{(\Phi_n)}(A\,;B)$, so by \cite{ERSZ2}, Lemma~4.2.III, there exists a $\psi \in D_0(h_{ \Phi_n \, 0})$
such that $0 \leq \psi \leq M$, $\psi|_B = 0$ and $\psi|_A \geq M - \varepsilon$.
Then $\chi \, \psi \in D(h_{\Phi_n \, 0}) \cap L_\infty$.
Moreover, it follows from Corollary~\ref{clbs106.5} and \cite{ERSZ2}, Lemma~3.3.II, that 
\begin{eqnarray*}
|||\widehat \ci^{(h_0)}_{\chi \, \psi}|||
 =  |||\widehat \ci^{(\Phi_n)}_{\chi \, \psi}||| 
& \leq & (1 + \varepsilon) \|\chi\|_\infty^2 \, |||\widehat \ci^{(\Phi_n)}_\psi|||
   + (1 + \varepsilon^{-1}) \|\psi\|_\infty^2 \, |||\widehat \ci^{(\Phi_n)}_\chi|||  \\[5pt]
& \leq & (1 + \varepsilon) + (1 + \varepsilon^{-1}) M^2 \, \lambda \, |||\widehat \ci^{(l)}_\chi|||  
\leq 1 + N \, \varepsilon
\end{eqnarray*}
where $N = 1 + 2 M^2 \lambda$ and $\lambda > 0$ is such that $h \leq l|_\gotA$.
So $(1 + N \, \varepsilon)^{-1/2} \chi \, \psi \in D_0(h_0)$.
Therefore 
\[
(1 + N \, \varepsilon)^{-1/2} (M - \varepsilon) 
\leq d_{(1 + N \, \varepsilon)^{-1/2} \chi \, \psi}(A\,;B)
\leq d^{(h_0)}(A\,;B)
\]
and $M \leq d^{(h_0)}(A\,;B)$.
Thus
\begin{equation}
d^{(h_0)}(A\,;B) 
= \inf_{n \in \Ni} d^{(\Phi_n)}(A\,;B)
= \lim_{n \to \infty} d^{(\Phi_n)}(A\,;B)
\;\;\; .
\label{ellbs112;1}
\end{equation}
If $d^{(h_0)}(A\,;B) = \infty$ then obviously $d^{(h_0)}(A\,;B) = d^{(\Phi_n)}(A\,;B)$ for all $n \in \Ni$.
Alternatively, if $d^{(h_0)}(A\,;B) < \infty$ then the lemma follows from Lemma~\ref{llbs111}
together with (\ref{ellbs112;1}).\hfill$\Box$

\ruimte

At this point we have control over the limits of the semigroups $S^{(\Phi)}$ and the distances
$d^{(\Phi)}$ associated with the truncated forms.

\ruimte

\noindent{\bf Proof of Theorem~\ref{tlbs121}\hspace{5pt} }\
Since, by hypothesis, $h\in\cc_l$ and  $l$ satisfies Condition~L  the foregoing results are applicable.
Let $\lambda > 0$ be such that $h \leq \lambda\,l|_\gotA$.

It follows from the Davies--Gaffney bounds Theorem~\ref{tdfd3.1} that 
\[
t \log (\one_A, S^{(h_0)}_t \one_B) \leq - 4^{-1} d^{(h_0)}(A\,;B)^2
\]
for all $t > 0$.
If $d^{(h_0)}(A\,;B) = \infty$ then (\ref{etlbs113;1}) is obviously valid.
So we may assume that $d^{(h_0)}(A\,;B) < \infty$.
Then in particular $A \neq \emptyset$.

Fix $x \in X$.
Let $\Phi_1,\Phi_2,\ldots \in \gotA$ be such that 
$0 \leq \Phi_1 \leq \Phi_2 \leq \ldots \leq 1$,  $\Phi_n|_{B^{(e)}(x;n)} = 1$
and $\supp \Phi_n$ is compact  for all $n \in \Ni$.
Let $\varepsilon > 0$.
By Lemma~\ref{llbs112} there exists an $n \in \Ni$ such that 
$d^{(h_0)}(A\,;B) = d^{(\Phi_n)}(A\,;B)$ and 
$r^2 = \lambda^{-1} \, d^{(e)}(A\,;B^{(e)}(x;n)^{\rm c})^2 > d^{(h_0)}(A\,;B)^2 + 8 \varepsilon$.
Then 
\[
|(\one_A, S^{(h_0)}_t \one_B) - (\one_A, S^{(\Phi_n)}_t \one_B)|
\leq 2\, e^{-(4t)^{-1} r^2} \, |A|^{1/2} \, |B|^{1/2}
\]
for all $t > 0$ by Corollary~\ref{clbs109}.
Moreover, by Corollary~\ref{clbs105.5}
there exists a $t_0 > 0$ such that 
\[
t \log (\one_A, S^{(\Phi_n)}_t \one_B)
\geq - 4^{-1} d^{(\Phi_n)}(A\,;B)^2 - \varepsilon
\]
for all $t \in \langle0,t_0]$.
Then
\begin{eqnarray*}
(\one_A, S^{(h_0)}_t \one_B)
& \geq & (\one_A, S^{(\Phi_n)}_t \one_B) - 2 e^{-(4t)^{-1} r^2} \, |A|^{1/2} \, |B|^{1/2}  \\[5pt]
& \geq & e^{-(4t)^{-1} d^{(\Phi_n)}(A;B)^2} \, e^{-\varepsilon t^{-1}} - 2 \,e^{-(4t)^{-1} r^2} \, |A|^{1/2} \, |B|^{1/2}  \\[5pt]
& \geq & e^{-(4t)^{-1} d^{(h_0)}(A;B)^2} \, e^{-\varepsilon t^{-1}} \Big(1 - 2\, e^{-\varepsilon t^{-1}} \, |A|^{1/2} \, |B|^{1/2} \Big)  
\end{eqnarray*}
for all $t \in \langle0,t_0]$, where we used  $d^{(h_0)}(A\,;B)^2 + 4 \varepsilon - r^2 \leq - 4 \varepsilon$
in the last step.
There is a $t_1 > 0$ such that 
$2 e^{-\varepsilon t_1^{-1}} \, |A|^{1/2} \, |B|^{1/2} \leq 4^{-1}$.
Then 
\[
(\one_A, S^{(h_0)}_t \one_B)
\geq 2^{-1} e^{-(4t)^{-1} d^{(h_0)}(A;B)^2} \, e^{-\varepsilon t^{-1}}
\]
and 
\[
t \log (\one_A, S^{(h_0)}_t \one_B)
\geq - 4^{-1} d^{(h_0)}(A\,;B)^2 - \varepsilon - t \log 2
\]
for all $t \in \langle0,t_0 \wedge t_1]$.
This completes the proof of  Theorem~\ref{tlbs121}.\hfill$\Box$

\ruimte

\noindent
{\bf Proof of Corollary~\ref{cas135}\hspace{5pt} }\
We may assume that $|A| \neq 0$ and $|B| \neq 0$.
If $\varphi,\psi \in L_2(X)$ then $P_A \psi \in L_2(A)$ and $P_B \varphi \in L_2(B)$.
So by the Davies--Gaffney bounds of Theorem~\ref{tdfd3.1} one deduces that 
\begin{eqnarray*}
|(\psi,P_A \, S_t \, P_B \varphi)|
& = & |(P_A \, \psi,S_t \, P_B \varphi)| \\[5pt]
& \leq & e^{-(4t)^{-1} \, d(A;B)^2} \, \|P_A \psi\|_2 \, \|P_B \varphi\|_2 
\leq e^{-(4t)^{-1} \, d(A;B)^2} \, \|\psi\|_2 \, \|\varphi\|_2 
\;\;\; .
\end{eqnarray*}
Hence 
\[
\limsup_{t \downarrow 0} t \log \|P_A S_t P_B\|_{2 \to 2} 
\leq - 4^{-1} d(A\,;B)^2  
\;\;\; .  \]
So it remains to show that 
\begin{equation}
- 4^{-1} d(A\,;B)^2  
\leq \liminf_{t \downarrow 0} t \log \|P_A S_t P_B\|_{2 \to 2} 
\label{ecas135;1}
\;\;\; .  
\end{equation}
Let $A_0 \subseteq A$ and $B_0 \subseteq B$ be measurable and relatively compact
and assume that $|A_0| > 0$ and $|B_0| > 0$.
Then 
\[
(\one_{A_0} , S_t \, \one_{B_0})
= (\one_{A_0} , P_A S_t P_B \, \one_{B_0})
\leq |A_0|^{1/2} \, \|P_A S_t P_B\|_{2 \to 2} \, |B_0|^{1/2}
\]
for all $t > 0$.
So
\begin{eqnarray*}
\liminf_{t \downarrow 0} t \log \|P_A S_t P_B\|_{2 \to 2} 
& \geq & \liminf_{t \downarrow 0} \Big( t \log (\one_{A_0} , S_t \, \one_{B_0}) - 2^{-1}t \log(|A_0| \, |B_0|) \Big)  \\[5pt]
& = & - 4^{-1} \, d(A_0\,;B_0)^2
\end{eqnarray*}
by Theorem~\ref{tlbs121}.
Now fix $x \in X$ and choose $A_0 = A \cap B^{(e)}(x\,;n)$ and 
$B_0 = B \cap B^{(e)}(x\,;n)$ with $n \in \Ni$ large enough.
Then
\begin{eqnarray*}
\liminf_{t \downarrow 0} t \log \|P_A S_t P_B\|_{2 \to 2} 
& \geq & - \lim_{n \to \infty} 4^{-1} \, d(A \cap B^{(e)}(x\,;n) \,; B \cap B^{(e)}(x\,;n))^2  \\[5pt]
& = & - 4^{-1} \, d(A\,;B)^2
\end{eqnarray*}
where the equality follows from Lemma~\ref{ldfd3n11}.
This proves (\ref{ecas135;1}) and the corollary.\hfill$\Box$

\ruimte

Now we turn to the proof of Theorem~\ref{tlbs122}.
First we use the continuity assumption  \ref{tlbs122}.\ref{tlbs122-1} of the theorem to identify
the distance.

\begin{lemma} \label{llbs141}
If $D_0(l) \subset C(X)$ then $d^{(l)}(A\,;B) = d^{(e)}(A\,;B)$
for all non-empty open subsets $A,B \subset X$.
\end{lemma}
\proof\
Let $M \in [0,d^{(l)}(A\,;B)] \cap \Ri$ and $\varepsilon > 0$.
Then by \cite{ERSZ2}, Lemma~4.2.III, there exists a $\psi \in D_0(l)$ such that 
$\psi|_B = 0$ almost everywhere and $\psi|_A \geq M - \varepsilon$ almost everywhere.
Since $D_0(l) \subset C(X)$ it follows that $\psi \in D(l)_{\loc} \cap C_{\rm b}(X)$.
Moreover, since $A$ and $B$ are open one has $\psi|_B = 0$ and $\psi|_A \geq M - \varepsilon$
pointwise.
So $|\psi(x) - \psi(y)| \geq M - \varepsilon$ for all $x \in A$ and $y \in B$.
It follows that $d^{(e)}(A\,;B) \geq M - \varepsilon$ and 
$d^{(e)}(A\,;B) \geq d^{(l)}(A\,;B)$.
The converse inequality is trivial (see Lemma~\ref{las352}.\ref{las352-1}).\hfill$\Box$

\ruimte

\noindent{\bf Proof of Theorem~\ref{tlbs122}\hspace{5pt} }\
Let $x_0,y_0 \in X$ and $\varepsilon \in \langle0,1]$.
Choose $A = B^{(e)}(x_0\,;\varepsilon)$ and $B = B^{(e)}(y_0\,;\varepsilon)$.
Let $t \in \langle0, 4^{-1} T]$.
Then by assumption \ref{tlbs122-2} of Theorem~\ref{tlbs122} and symmetry of the kernel
\begin{eqnarray*}
K_t(x_0\,;y_0)
& \leq & K_{(1+\varepsilon) t}(x_0\,;y) \, (1+\varepsilon)^\nu \, e^{\omega(1 + \varepsilon t^{-1})}  \\[5pt]
& \leq & K_{(1+\varepsilon)^2 t}(x\,;y) \, (1+\varepsilon)^{2\nu} \, e^{2\omega(1 + \varepsilon t^{-1})}  
\end{eqnarray*}
for all $x \in A$ and $y \in B$.
So
\begin{eqnarray*}
K_t(x_0\,;y_0)
& \leq & |A|^{-1} \, |B|^{-1} \int_A dx \int_B dy \, 
    K_{(1+\varepsilon)^2 t}(x\,;y) \, (1+\varepsilon)^{2\nu} \, e^{2\omega(1 + \varepsilon t^{-1})}  \\[5pt]
& = & |A|^{-1} \, |B|^{-1} \, (1+\varepsilon)^{2\nu} \, e^{2\omega(1 + \varepsilon t^{-1})} \,
        (\one_A \, S^{(l)}_{(1+\varepsilon)^2 t} \one_B)
\;\;\; .
\end{eqnarray*}
Therefore
\[
t \log K_t(x_0\,;y_0)
\leq t \log (\one_A \, S^{(l)}_{(1+\varepsilon)^2 t} \one_B)
   - t \log(|A| \, |B|)
   + 2\,\nu \, t \log(1+\varepsilon)
   + 2\, \omega \, t 
   + 2 \,\omega \, \varepsilon
\]
for all $t \in \langle0, 4^{-1} T]$.
Since $|A|, |B| > 0$ one deduces from Corollary~\ref{cas2.1} that 
\[
\limsup_{t \downarrow 0} \,t \log K_t(x_0\,;y_0)
\leq - 4^{-1} (1+\varepsilon)^{-2} d^{(l)}(A\,;B) + 2\, \omega \, \varepsilon
= - 4^{-1} (1+\varepsilon)^{-2} d^{(e)}(A\,;B) + 2 \,\omega \, \varepsilon
\]
where we used Lemma~\ref{llbs141} in the last step.
But 
\[
d^{(e)}(x_0,y_0) - 2 \,\varepsilon
\leq d^{(e)}(B^{(e)}(x_0\,;\varepsilon) \,; B^{(e)}(y_0\,;\varepsilon))
\leq d^{(e)}(x_0,y_0) + 2\, \varepsilon
\]
by the triangle inequality.
So taking the limit $\varepsilon \downarrow 0$ one establishes that 
\[
\limsup_{t \downarrow 0} \,t \log K_t(x_0\,;y_0)
\leq - 4^{-1} d^{(e)}(x_0,y_0)
\;\;\; .  \]
Since, by an analogous argument,
\[
K_t(x_0\,;y_0)
\geq K_{(1-\varepsilon)^2 t}(x\,;y) \, (1+\varepsilon)^{-2\nu} \, e^{-3\omega(1 + \varepsilon t^{-1})}  
\]
for all $t \in \langle 0,T]$, $\varepsilon \in \langle0,1/2\rangle$,
$x \in B^{(e)}(x_0\,;\varepsilon)$ and $y \in B^{(e)}(y_0\,;\varepsilon)$
one deduces similarly that 
\[
\liminf_{t \downarrow 0} \,t \log K_t(x_0\,;y_0)
\geq - 4^{-1} d^{(e)}(x_0,y_0)
\;\;\; .  \]
So $\lim_{t \downarrow 0}\, t \log K_t(x_0\,;y_0) = - 4^{-1} d^{(e)}(x_0,y_0)$
and Theorem~\ref{tlbs122} follows.\hfill$\Box$

\section{Applications}\label{Sas4}

 Theorems~\ref{tlbs121} and \ref{tlbs122} have a broad range of applications to 
second-order,  divergence-form, elliptic operators both degenerate and non-degenerate.
First we discuss the application to operators on $\Ri^d$.

Let $c_{ij}$ be bounded real-valued measurable functions on $\Ri^d$
and assume  that the 
$d\times d$-matrix $C=(c_{ij})$ is symmetric and positive-definite almost-everywhere.
Define  $h$ by
\[
h(\varphi)=\sum^d_{i,j=1}(\partial_i\varphi,c_{ij}\partial_j\varphi)
\]
where $\partial_i=\partial/\partial x_i$ and $D(h)=W^{1,2}(\Ri^d)$.
Then  $h$ is a positive form which is strongly local and Markovian.
It is not in general closed or even closable.
A characterization of closable forms in one dimension can be found in \cite{FOT}, 
Theorem~3.1.6, and sufficient conditions  in higher  dimensions are given in 
 \cite{FOT}, Section~3.1 (see also \cite{RW} or \cite{MR}, Chapter~II).
Nevertheless, $h\leq \lambda\,l$ 
where  $\lambda$ denotes the essential supremum of the matrix norms $\|C(x)\|$
and $l$ is the form of the usual Laplacian on $\Ri^d$, i.e.,
$D(l)=W^{1,2}(\Ri^d)$ and 
\[
l(\varphi)=\sum^d_{i=1}\|\partial_i\varphi\|_2^2=\|\nabla\varphi\|_2^2
\]
for all $\varphi\in D(l)$.

First, consider the case of strongly elliptic operators, i.e., assume that
 $C\geq\mu\, I$ almost everywhere, with $\mu>0$.
Then one has $\lambda\,l\geq h\geq\mu\,l$ and the form $h$ is closed on $W^{1,2}(\Ri^d)$.
It follows readily that $h$ is a diffusion satisfying Condition~$L$ and the distance $d^{(e)}$
 is the Riemannian distance corresponding to the metric $C^{-1}$.
Moreover, Conditions~\ref{tlbs122-1} and~\ref{tlbs122-2} of Theorem~\ref{tlbs122}
are valid. 
The latter condition follows from the parabolic Harnack inequality (see, for example,
\cite{Sal5}, Chapter~5 and in particular Corollary~5.4.6).
Therefore Varadhan's small time asymptotic result (\ref{eas1.1}) follows for each strongly
elliptic operator by Corollary~\ref{cas2.1} and Theorem~\ref{tlbs122}.

Secondly,  consider degenerate operators.
Thus $C\geq0$ almost everywhere but this is the only coercivity condition.
Nevertheless, $h$ is a strongly local Markovian form and  $h\leq\lambda\,l|_\gotA$.
Therefore Theorem~\ref{tlbs121} is directly applicable to $h$ with no further assumptions
on the coefficients.
It then follows that $h_0$ is strongly local and the asymptotic identification
(\ref{etlbs113;1}) is valid.
Further detail on the asymptotic behaviour requires more detailed analysis of the 
set-theoretic distance $d^{(h_0)}(A\,;B)$.
In the degenerate situation the
set $U_h=\{x\in\Ri^d:\,\|C(x)\|=0\}$ may have non-zero measure and the diffusion
associated with the relaxation $h_0$ occurs on the closure of the  complement 
$\Ri^d\backslash U_h$.
In fact the effective space of the diffusion can be smaller 
since sets of capacity zero act as obstacles \cite{RSi}.
\smallskip

The properties of the distance $d^{(h)}(A\,;B)$ can be understood  for non-degenerate,
or weakly degenerate, elliptic operators.
In particular it can be analyzed for subelliptic operators in quite general situations.

Let $M$ be a Riemannian manifold  and $X_1,\ldots, X_n$ 
smooth vector fields on $M$.
Each such vector field $X$ defines a closed linear partial differential operator,
also denoted by $X$, on $L_2(M)$.
Now consider the form 
\[
h(\varphi)=\sum^d_{i=1}\|X_i\varphi\|_2^2=\| \, |\nabla\varphi| \, \|_2^2
\]
with $D(h)=\bigcap^n_{i=1}D(X_i)$ and  $\nabla\varphi=(X_1\varphi,\ldots,X_n\varphi)$.
The form is automatically closed regular and strongly local.
Next  assume the form is densely-defined and that the vector fields satisfy
the H\"ormander condition of order $r$.
It follows readily that
$\psi\in D_0(h)$ if and only if $|\nabla\psi|\in L_\infty(M)$ and 
$|\nabla\psi|^2\leq 1$ almost everywhere.
But it then follows from \cite{RS}, Theorem~17, that $\psi$ is locally Lipschitz.
In fact it has a Lipschitz derivative of order $1/r$.
Consequently, $D_0(l)\subseteq C(M)$.
Thus the first assumption of Theorem~\ref{tlbs122} is verified in quite general circumstances.
The second assumption is also verifiable in many situations as it is a consequence of a 
parabolic Harnack inequality.
As a specific illustration we give an application to subelliptic operators on Lie groups.

\begin{prop} \label{pas461}
Let $X_1,\ldots,X_n$ be right invariant vector fields on a Lie group 
$G$ which satisfy the H\"ormander condition, i.e., the vector fields
generate the Lie algebra of $G$.
For all $i,j \in \{ 1,\ldots,n \} $ let $c_{ij} = c_{ji} \in L_\infty(G)$ be 
real valued and assume there is a $\mu > 0$ such that
$(c_{ij}(g)) \geq \mu \, I$ for almost every $g \in G$.
Define the quadratic form $h$ on $L_2(G)$ by
\[
h(\varphi) 
= \sum_{i,j=1}^n   (X_i \varphi,  c_{ij}\,X_j \varphi)
\]
with  form domain $D(h) = \bigcap_{i=1}^n D(X_i)$.
Then the form $h$ is a diffusion
satisfying {\rm Condition~L} and the assumptions 
of Theorem~{\rm \ref{tlbs122}} are satisfied.
Hence 
\[
\lim_{t \downarrow 0} \,t \log K_t(g_1\,;g_2) = - 4^{-1} d(g_1\,;g_2)^2
\]
for all $g_1,g_2 \in G$, where $K$ is the kernel of the semigroup
generated by the operator $H$ associated to the form $h$ and 
\begin{equation}
d(g_1\,;g_2)
= \sup \{ |\psi(g_1) - \psi(g_2)| : \psi \in D(h)_{\rm loc} \cap C_{\rm b}(G) 
         \mbox{ and } |||\widehat \ci_\psi||| \leq 1 \}
\;\;\; .  
\label{epas461;1}
\end{equation}
\end{prop}
\proof\
It is straightforward to check that $h$ is a diffusion.
Note also that $h$ is bounded above and below by multiples of  the form
\[
l(\varphi)=\sum^n_{i=1}\|X_i\varphi\|_2^2
\]
associated with the sublaplacian given by the vector fields.
Therefore the distances $d^{(e)}$ for $l$ and $d$ (defined by (\ref{epas461;1}))
for $h$ are equivalent.
Condition~L and the first assumption of Theorem~\ref{tlbs122} 
for the form $l$ and then also for the form $h$
follow since the vector 
fields satisfy the H\"ormander condition.
The second assumption of the theorem follows from the parabolic Harnack inequality 
(again see \cite{Sal5}, Chapter~5 and in particular Corollary~5.4.6,
or the first remark following Theorem~4.4 on page~35--36 in \cite{Sal2}).\hfill$\Box$

\ruimte

Next we note that Norris' result on the small time behaviour of the Laplace--Beltrami
operator can also be deduced from our results.
Let $M$ be a $d$-dimensional Lipschitz Riemannian manifold with Borel measure $\mu$ which 
 under some chart is  locally equivalent with the Lebesgue measure.
Further let $l$ denote the Dirichlet form
\[
l(\varphi)=\int_Md\mu\,|\nabla\varphi|^2
\]
where $\nabla$ denotes the usual gradient and $D(l)=W^{1,2}(M)$.
The corresponding self-adjoint operator on $L_2(M)$ is the 
 Laplace--Beltrami operator for $M$.
Using a local coordinate chart one easily argues that  
 $\psi\in W^{1,2}_{\rm loc}\cap L_\infty$ and $|||\widehat \ci^{(l)}_\psi|||\leq 1$ 
if and only if $\psi\in W^{1,\infty}$ and 
$\|\nabla\psi\|_\infty\leq 1$.
In particular $D_0(l)\subseteq C(M)$.
Since the Riemannian distance on $M$ is given by
\[
d(x\,;y)=\sup\{|\psi(x)-\psi(y)|: \psi\in W^{1,\infty}\,,\,\|\nabla \psi\|_\infty\leq 1\}
\]
it follows that Condition~L is satisfied.
Moreover, the parabolic Harnack inequality is valid (and is used in the proof of Norris
(see \cite{Nor1}, page~87)).
Therefore, by the above reasoning, one can apply Theorem~\ref{tlbs122} 
to establish 
\[
\lim_{t \downarrow 0}\,t\log K_t(x\,;y)=- 4^{-1} d(x\,;y)^2
\]
for all $x,y \in M$.
Note that as in Norris' argument this proof requires neither smoothness 
nor completeness of $M$.

\smallskip

Finally we note that the asymptotic estimates can be extended to operators
with lower order terms by various arguments such as perturbation theory.
For example, if  $h$ is a Dirichlet form and $v$  the 
form of a real, positive,  bounded multiplication operator~$V$
then the Trotter product formula
\[
S^{(h+v)}_t=\lim_{n\to\infty}(S^{(h)}_{t/n}S^{(v)}_{t/n})^n
\;\;\;.
\]
Therefore one deduces from positivity of the semigroups that
\begin{equation}
S^{(h)}_t\,e^{-t\|V\|_\infty}\leq S^{(h+v)}_t
\leq S^{(h)}_t
\label{eas4.444}
\end{equation}
for all $t > 0$.
Hence
\begin{equation}
\lim_{t \downarrow 0}\,t\log(\one_A,S^{(h+v)}_t\one_B)=\lim_{t \downarrow 0}\,t\log(\one_A,S^{(h)}_t\one_B)
\label{eas4;76}
\end{equation}
for all relatively compact  measurable subsets $A$ and $B$.
Thus  the small time asymptotics is independent of $v$.

Next, if $h$ is a regular Dirichlet for and $v$ is the form of the operator of multiplication by a real, positive,
locally bounded measurable function $V$, then obviously $h+v$ is densily defined.
But in  addition $h+v$ is closed by \cite{bSim5} Theorem~4.1, since the closed forms
$h+v_n$ converge monotonically upwards to $h+v$, where 
$v_n$ denotes the form of the  operator of multiplication by the bounded function $V \, \one_{X_n}$
and $X_1 \subset X_2 \subset \ldots$ are measurable subsets of $X$ with $X = \bigcup_{n=1}^\infty X_n$.

These observations extend to the following result.

\begin{prop}\label{pas471}
Let $h$ be a positive form satisfying the hypotheses of Theorem~$\ref{tlbs121}$
and  with $D(h) = D(l)$.
Further, let  $v$  be the form of the operator of multiplication by a real, positive,
locally bounded measurable function $V$.
Then
\[
(h+v)_0=h_0+v=\rlim_{\varepsilon\to0} \,(h_\varepsilon+v)
\]
and
\[
\lim_{t \downarrow 0}\,t\log(\one_A,S^{((h+v)_0)}_t\one_B)=-4^{-1}d^{(h_0)}(A\,;B)^2
\]
for all relatively compact  measurable subsets $A$ and $B$.
\end{prop}
\proof\
First, $h_0+v$ is closed, densely-defined and $h_0+v\leq h+v\leq h_\varepsilon+v$ for all 
$\varepsilon > 0$, where in the second inequality we use the additional assumption $D(h) = D(l)$.
(Recall that by definition $D(h_\varepsilon) = D(l)$.)
Therefore $h_0+v\leq (h+v)_0\leq \hat h$ where $\hat h=\rlim_{\varepsilon \downarrow 0} (h_\varepsilon+v)$.
We now establish the first statement of the proposition by proving that $\hat h=h_0+v$.

Let $\Omega_1 \subset \Omega_2 \subset \ldots$ be open relatively compact subsets of $X$ 
such that $X = \bigcup_{n=1}^\infty \Omega_n$.
For all $n \in \Ni$ 
let $v_n$ denote the form of the  operator of multiplication by the bounded function $V \, \one_{\Omega_n}$.
Let $\varepsilon > 0$.
Set $h_{\varepsilon,n}=h_\varepsilon+v_n$ for all $n \in \Ni$.
Then $S^{(h_{\varepsilon,n})}_t\leq S^{(h_{\varepsilon})}_t$ for all $t>0$ by (\ref{eas4.444}) and 
\[
(\one_A,S^{(h_{\varepsilon,n})}_t\one_B)\leq (\one_A,S^{(h_{\varepsilon})}_t\one_B)\leq
 e^{-(\lambda+\varepsilon)^{-1/2} d^{(e)}(A;B)^2(4t)^{-1}}|A|^{1/2}|B|^{1/2}
\]
for all non-empty open $A,B \subset X$ and $t > 0$ where the second bound follows by 
Lemmas~\ref{llbs232} and \ref{las352}.\ref{las352-1} and $\lambda > 0$ is such that $h \leq \lambda \, l|_\gotA$.
Hence the  positive self-adjoint operator $H_{\varepsilon,n}$
associated with $h_{\varepsilon,n}$ has a finite speed of propagation.
Explicitly
\[
(\psi, \cos(t H_{\varepsilon,n}^{1/2}) \varphi) = 0
\]
for all non-empty open $A,B \subset X$, all $\varphi \in L_2(A)$, $\psi \in L_2(B)$
and  all $t \in \Ri$ with $|t| \leq (\lambda+\varepsilon)^{-1/2} \, d^{(e)}(A\,;B)$ by \cite{ERSZ1}, Lemma~3.3,
and the assumptions on $h$. 
Then arguing as in the proof of Proposition~\ref{plbs106}
one deduces that 
\[
\cos(t H_{\varepsilon,n}^{1/2}) \varphi
= \cos(t H_{\varepsilon,m}^{1/2}) \varphi
\]
for all $n,m \in \Ni$, measurable $A \subset X$, $\varphi \in L_2(A)$ and $t \in \Ri$ with 
$m \geq n$, $\emptyset \neq \overline A \subset \Omega_n$ and
$|t| \leq (\lambda+\varepsilon)^{-1/2} d^{(e)}(A\,;\Omega_n^{\rm c})$.
The proof is a repetition of the arguments used to prove Proposition~\ref{plbs106} but now
one also uses locality of $v_n$ and $v_m$.
Therefore using the representation (\ref{eas3.297}) one obtains an estimate
\[
\|S^{(h_{\varepsilon,n})}_t \varphi - S^{(h_{\varepsilon,m})}_t \varphi\|_2
\leq 2 \,e^{-(4t)^{-1} r_n^2} \, \|\varphi\|_2
\]
for all $n,m \in \Ni$, measurable $A \subset X$, $\varphi \in L_2(A)$ and $t>0$ with 
$m \geq n$, $\emptyset \neq \overline A \subset \Omega_n$, where
$r_n = (\lambda+\varepsilon)^{-1/2} \, d^{(e)}(A\,; \Omega_n^{\rm c})$.
Now if $m\to\infty$ then $h_{\varepsilon,m}$ converges monotonically upward to $h_\varepsilon+v$.
Hence  $S^{(h_{\varepsilon,m})}_t$ converges strongly to $S^{(h_{\varepsilon}+v)}_t$ by \cite{bSim5}, Theorem~3.1.
But  $S^{(h_{\varepsilon}+v)}_t$ converges strongly to $S^{(\hat h)}_t$ as  $\varepsilon \downarrow 0$
and  $S^{(h_{\varepsilon,n})}_t$ converges strongly
to $S^{(h_0+v_n)}_t$ as $\varepsilon \downarrow 0$ because $\rlim_{\varepsilon\downarrow 0} h_{\varepsilon,n} = h_0+v_n$.
Combining these observations one deduces that 
\begin{equation}
\|S^{(h_0+v_n)}_t \varphi - S^{(\hat h)}_t \varphi\|_2
\leq 2 \,e^{-(4t)^{-1} r_n^2} \, \|\varphi\|_2
\label{eas4;75}
\end{equation}
for all $n \in \Ni$, measurable $A \subset X$, $\varphi \in L_2(A)$ and $t>0$
with $\emptyset \neq \overline A \subset \Omega_n$.
Finally $\rlim_{n\to\infty} (h_0+v_n)= h_0+v$ by monotone convergence and as $r_n\to\infty$ as $n\to\infty$ 
 one concludes that $S^{(h_0+v)}_t \varphi =S^{(\hat h)}_t \varphi$ for all 
$\varphi \in L_{2,c}(X)$ and $t > 0$.
Thus  $S^{(h_0+v)}_t  =S^{(\hat h)}_t $ for  all  $t > 0$.
Hence $\hat h=h_0+v$.

Finally we deduce from (\ref{eas4.444})  that 
\[
(\one_A,S^{(\hat h)}_t\one_B)
= (\one_A,S^{(h_0+v)}_t\one_B)
\leq (\one_A,S^{(h_0)}_t\one_B)
\leq e^{-d^{(h_0)}(A;B)^2(4t)^{-1}}|A|^{1/2}|B|^{1/2}
\]
for all non-empty open $A,B \subset X$ and $t > 0$, where the last bound uses Theorem~\ref{tdfd3.1}.
Finally the proof of the proposition is completed by repeating the arguments used above
in the proof of Theorem~\ref{tlbs121} with $S^{(h_0)}$ replaced by $S^{(\hat h)}$ and 
$S^{(\Phi)}$ replaced by $S^{(h_0+v_n)}$.
The equality (\ref{eas4;76}) and inequality (\ref{eas4;75}) replace Corollaries~\ref{clbs105.5} and 
\ref{clbs109}.\hfill$\Box$

\ruimte

Note that the identification $(h+v)_0=h_0+v$ shows that $(h+v)_0$ is local in the sense of \cite{FOT}.

\section*{Acknowledgements}

This work  was  supported by the Australian Research Council (ARC)
Discovery Grant DP 0451016.
The greater part of this work was carried out during a visit of the first named
author to the Australian National University.
The work was completed whilst the first and second named authors were visiting the
Centre International de Rencontres Math\'ematiques at the Universit\'e de la M\'editerran\'ee
at Luminy.

\end{document}